\documentclass[12pt]{article}
\usepackage{latexsym,amssymb,amsmath}
\begin{document}

\baselineskip=20pt

\newcommand{\rd}{\mbox{Rad}}
\newcommand{\kn}{\mbox{ker}}
\newcommand{\psp}{\vspace{0.4cm}}
\newcommand{\pse}{\vspace{0.2cm}}
\newcommand{\ptl}{\partial}
\newcommand{\dlt}{\delta}
\newcommand{\sgm}{\sigma}
\newcommand{\al}{\alpha}
\newcommand{\be}{\beta}
\newcommand{\G}{\Gamma}
\newcommand{\gm}{\gamma}
\newcommand{\lmd}{\lambda}
\newcommand{\td}{\tilde}
\newcommand{\vf}{\varphi}
\newcommand{\ad}{\mbox{ad}}
\newcommand{\stl}{\stackrel}
\newcommand{\ol}{\overline}
\newcommand{\es}{\epsilon}
\newcommand{\ves}{\varepsilon}
\newcommand{\la}{\langle}
\newcommand{\ra}{\rangle}

\begin{center}{\Large \bf Generalizations of the Block Algebras}
\end{center}
\vspace{0.2cm}

\begin{center}{\large Xiaoping Xu}\end{center}
\begin{center}{Department of Mathematics, The Hong Kong University of Science \& Technology}\end{center}
\begin{center}{Clear Water Bay, Kowloon, Hong Kong}\footnote{Research supported
 by Hong Kong Research Grants Council Competitive Earmarked Research Grant HKUST709/96P.}\end{center}

\vspace{0.3cm}

\begin{center}{\Large \bf Abstract}\end{center}
\vspace{0.2cm}

{\small In this paper, we introduce two new families of infinite-dimensional simple Lie algebras and a new family of infinite-dimensional simple Lie superalgebras. These algebras can be viewed as generalizations of the Block algebras.}

\section{Introduction}

 Block [B] introduced certain analogues of the Zassenhaus algebras over a field with characteristic 0  and showed that they are simple under a certain condition. A family of infinite-dimensional simple Lie superalgebras was constructed in [X2] as the Lie superalgebras induced by the Hamiltonian superoperators associated with certain Novikov-Poisson algebras (cf. [X1]). These algebras can be viewed as generalizations of the Block algebras. In the classification of Gel'fand-Dorfman bialgebras (equivalently, certain Hamiltonian pairs) whose Novikov algebraic structures are simple with an idempotent element, we found in [X3] new Lie algebras, which can also be viewed as  generalizations of the Block algebras . In this paper, we shall present more general simple Lie superalgebras that can be viewed as generalizations of the Block algebras. One of our motivations of this paper is to understand the simplicity of the Lie superalgebras generated by certain quadratic conformal superalgebras (cf. [X3], [X5]). This is a subsequent paper to [X4]. 

Throughout this paper, let $\Bbb{F}$ be a field  with characteristic 0. All the vector spaces are assumed over $\Bbb{F}$. Denote by $\Bbb{Z}$ the ring of integers, by $\Bbb{N}$ the set of natural numbers $\{0,1,2,3,...\}$.

Let $\G$ be a nonzero torsion-free abelian group. Moreover, we also view $\G$ as a $\Bbb{Z}$-module. Suppose that $\phi(\cdot,\cdot):\G\times \G\rightarrow \Bbb{F}$ is a skew-symmetric $\Bbb{Z}$-bilinear form and $\vf:\G\rightarrow \Bbb{F}$ is an additive group homomorphism. Let ${\cal A}_{\G}$ be a vector space with a basis $\{x^{\al}\mid \al\in\G\}$. The following Lie bracket on ${\cal A}_{\G}$ is given by Block [B]:
$$[x^{\al},x^{\be}]_B=(\phi(\al,\be)+\vf(\al-\be))x^{\al+\be}\qquad\mbox{for}\;\;\al,\be\in\G,\eqno(1.1)$$
where either $\vf\equiv 0$ or there exists another additive group homomorphism $\vf_1:\G\rightarrow \Bbb{F}$ such that 
$$\phi(\al,\be)=\vf (\al)\vf_1(\be)-\vf(\be)\vf_1(\al)\qquad\mbox{for}\;\;\al,\be\in\G.\eqno(1.2)$$
We call $({\cal A}_{\G},[\cdot,\cdot]_B)$ a {\it Block algebra}. If $\vf\equiv 0$ and $\mbox{Rad}_{\phi}=\{0\}$, then $({\cal A}_{\G},[\cdot,\cdot]_B)$ has the one-dimensional center $\Bbb{F}x^0$, and ${\cal A}_{\G}/\Bbb{F}x^0$ is simple. When $\phi$ is given by (1.2),
$$\kn_{\vf}\bigcap\kn_{\vf_1}=\{0\}\;\;\mbox{and}\;\;\vf_1(\al)\neq 2\;\;\mbox{for}\;\;\al\in\kn_{\vf},\eqno(1.3)$$
the Lie algebra $({\cal A}_{\G},[\cdot,\cdot]_B)$ is simple (cf. [B]). The case when the second condition in (1.3) fails was tackled by Dokovic and Zhao [DZ]. Some special cases of our generalized Lie algebras of Cartan type H in [X4] can be viewed as generalizations  of $({\cal A}_{\G},[\cdot,\cdot]_B)$ with $\vf\equiv 0$. 
In this paper, we shall introduce three families of generalizations of the Block algebra when $\phi$ is given by (1.2). 

Let $n$ be a positive integer. We denote 
$$\ol{1,n}=\{1,2,3,...,n\}.\eqno(1.4)$$
 Pick
$${\cal J}_p\in\{\{0\},\Bbb{N}\}\qquad\mbox{for}\;\;p\in\ol{1,n}.\eqno(1.5)$$
Let $\{\vf_p\mid p\in\ol{1,n}\}$ be $n$ additive group homomorphisms from $\G$ to $\Bbb{F}$ such that
$$\vf_p\not\equiv 0\;\;\mbox{if}\;\;{\cal J}_p=\{0\}\qquad\mbox{for}\;\;p\in\ol{1,n},\eqno(1.6)$$
$$\bigcap_{p=1}^n\kn_{\vf_p}=\{0\}.\eqno(1.7)$$
Condition (1.7) implies that $\G$ is isomorphic to an additive subgroup of $\Bbb{F}^n$.

Set 
$$\vec{\cal J}={\cal J}_1\times{\cal J}_2\times\cdots\times{\cal J}_n,\eqno(1.8)$$
where the addition on $\vec{\cal J}$ is defined componentwisely. Moreover, we denote: 
$$i_{[p]}=(0,..,\stackrel{p}{i},0,...,0)\qquad\mbox{for}\;\;i\in{\cal J}_p.\eqno(1.9)$$
We assume that $\G$ is a  torsion-free abelian group throughout this paper. We allow $\G=\{0\}$. Let ${\cal A}_n$ be a vector space  with a basis 
$$\{x^{\al,\vec{i}}\mid (\al,\vec{i})\in\G\times\vec{\cal J}\}.\eqno(1.10)$$
We define an algebraic operation $\cdot$ on ${\cal A}_n$ by:
$$x^{\al,\vec{i}}\cdot x^{\be,\vec{j}}=x^{\al+\be,\vec{i}+\vec{j}}\qquad\mbox{for}\;\;(\al,\vec{i}),(\be,\vec{j})\in\G\times\vec{\cal J}.\eqno(1.11)$$
Then ${\cal A}_n$ forms a commutative associative algebra with an identity element $x^{0,\vec{0}}$, which will be simply denoted by $1$ in the rest of this paper. Throughout this paper, the notion ``$\cdot$'' will be invisible in a product when the context is clear. For $p\in\ol{1,n}$, we define $\ptl_p\in \mbox{End}\:{\cal A}_n$ by
$$\ptl_p(x^{\al,\vec{i}})=\vf_p(\al)x^{\al,\vec{i}}+i_px^{\al,\vec{i}-1_{[p]}}\eqno(1.12)$$
for $(\al,\vec{i})\in\G\times\vec{\cal J}$, where we adopt the convention that if a notion is not defined but technically appears in an expression, we always treat it as zero; for instance, $x^{\al,-1_{[1]}}=0$ for any $\al\in\G$.  It can be verified that $\{\ptl_p\mid p\in\ol{1,n}\}$ are derivations of ${\cal A}_n$. 
\psp

{\it Class I}. 
\psp

Let $n=2$. We define the following algebraic operation $[\cdot,\cdot]_1$ on ${\cal A}_2$:
$$[u,v]_1=\ptl_1(u)\ptl_2(v)-\ptl_1(v)\ptl_2(u)+u\ptl_1(v)-\ptl_1(u)v\qquad\mbox{for}\;\;u,v\in{\cal A}_2.\eqno(1.13)$$
Then  $({\cal A}_2,[\cdot,\cdot]_1)$ forms a Lie algebra by (3.45) in [X3]. By (1.7), there exist at most one element $\sgm_1\in\kn_{\vf_1}$ and one element $\sgm_2\in\kn_{\vf_1}$ such that $\vf_2(\sgm_1)=1$ and $\vf_2(\sgm_2)=2$, respectively. If such a $\sgm_1$ exists, then $x^{\sgm_1,\vec{0}}$ is a central element of ${\cal A}_2$.  We treat $x^{\sgm_i,\vec{0}}$ as $0\in {\cal A}_2$ when such  $\sgm_i$ do not exist. Form a quotient Lie algebra
$${\cal B}_2={\cal A}_2/\Bbb{F}x^{\sgm_1,\vec{0}}\eqno(1.14)$$
whose induced Lie bracket is still denoted by $[\cdot,\cdot]_1$ when the context is clear. 

\psp

{\bf Theorem 1}. {\it The Lie algebra} $({\cal B}_2,[\cdot,\cdot]_1)$ {\it is simple if} $\vec{\cal J}\neq\{\vec{0}\}$ {\it or} $\sgm_2$ {\it does not exists}. If $\vec{\cal J}=\{\vec{0}\}$ and $\sgm_2$ {\it exists, then} ${\cal B}_2^{(1)}=[{\cal B}_2,{\cal B}_2]$ {\it is simple and} ${\cal B}_2={\cal B}_2^{(1)}\oplus (\Bbb{F}x^{\sgm_2,\vec{0}}+\Bbb{F}x^{\sgm_1,\vec{0}})$.
\psp

The special case $\vec{\cal J}=\{\vec{0}\}$ of Theorem 1 is due to Block [B] and due to Dokovic and  Zhao [DZ]. However, our proof is different from theirs. We also remark that if we allow $\ptl_2=0$ (equivalently, (1.6) only holds for $p=1$), then $({\cal A}_2,[\cdot,\cdot])$ is a known rank-one simple Lie algebra of generalized Witt type, which including the centerless Virasoro algebras as a special case. 
\psp

{\it Class II}.
\psp

Let $n=4$ and assume
$$\bigcap_{p\neq q\in\ol{1,4}}\kn_{\vf_q}\setminus\kn_{\vf_p}\neq\emptyset\qquad\mbox{for}\;\;p\in\ol{1,4}.\eqno(1.15)$$
Let
$$\al_0\in\G\setminus[(\kn_{\vf_1}\bigcap\kn_{\vf_2})\bigcup(\kn_{\vf_3}\bigcap\kn_{\vf_4})]\eqno(1.16)$$
 be a fixed  element such that
$$\vf_p(\al_0)\in\vf_p(\kn_{\vf_3}\bigcap\kn_{\vf_4}),\;\;\vf_q(\al_0)\in\vf_q(\kn_{\vf_1}\bigcap\kn_{\vf_2})\eqno(1.17)$$
for $p=1,2$ and $q=3,4$. We define the following algebraic operation $[\cdot,\cdot]_2$ on ${\cal A}_4$ by
\begin{eqnarray*}[u,v]_2&=&x^{\al_0,\vec{0}}(\ptl_1(u)\ptl_2(v)-\ptl_1(v)\ptl_2(u))+(\ptl_3(u)\ptl_4(v)-\ptl_3(v)\ptl_4(u))\\& &+\vf_3(\al_0)(u\ptl_4(v)-\ptl_4(u)v)-\vf_4(\al_0)(u\ptl_3(v)-\ptl_3(u)v)\hspace{3.4cm}(1.18)\end{eqnarray*}
for $u,v\in{\cal A}_4$. It can be verified that $({\cal A}_4,[\cdot,\cdot]_2)$ forms a Lie algebra.
By (1.7), there exists at most one element $\sgm\in\kn_{\vf_1}\bigcap\kn_{\vf_2}$ such that $\vf_3(\sgm+\al_0)=\vf_4(\sgm+\al_0)=0$. If such a $\sgm$ exists, then $x^{\sgm,\vec{0}}$ is a central element of ${\cal A}_4$. Form a quotient Lie algebra
$${\cal B}_4={\cal A}_4/\Bbb{F}x^{\sgm,\vec{0}}\eqno(1.19)$$
whose induced Lie bracket is still denoted by $[\cdot,\cdot]_2$ when the context is clear, where we treat $x^{\sgm,\vec{0}}$ as $0\in {\cal A}_4$ when such a $\sgm$ does not exist. 
\psp

{\bf Theorem 2}. {\it The Lie algebra} $({\cal B}_4,[\cdot,\cdot]_2)$ {\it is simple if} $\vec{\cal J}\neq \{\vec{0}\}$ {\it or there does not exist} $\rho\in\G$ {\it such that} 
$$\vf_1(\rho-\al_0)=\vf_2(\rho-\al_0)=\vf_3(\rho+2\al_0)=\vf_3(\rho+2\al_0)=0.\eqno(1.20)$$
{\it If} $\vec{\cal J}=\{\vec{0}\}$ {\it and there exists} $\rho\in\G$ {\it such that (1.20) holds (it is unique by (1.7)), then } ${\cal B}_4^{(1)}=[{\cal B}_4,{\cal B}_4]_2$ {\it is simple and} ${\cal B}_4={\cal B}_4^{(1)}\oplus (\Bbb{F}x^{\rho,\vec{0}}+\Bbb{F}x^{\sgm,\vec{0}})$. 

\psp

{\it Class III}.
\psp

The third class of generalized Block algebras are Lie superalgebras, where some special cases of these Lie superalgebras had been obtained in [X2]. 
Set
$$\td{\cal A}={\cal A}_2\times {\cal A}_2=\td{\cal A}_0\oplus \td{\cal A}_1\eqno(1.21)$$
with
$$\td{\cal A}_0=({\cal A}_2,0),\qquad \td{\cal A}_1=(0,{\cal A}_2).\eqno(1.22)$$
Moreover, we denote
$$u_{[0]}=(u,0),\;\;u_{[1]}=(0,u)\eqno(1.23)$$
for $u\in {\cal A}_2.$

Let $\al_0\in \G$ and $\ves\in \Bbb{F}$ be  fixed elements.
 Define an algebraic operation $[\cdot,\cdot]$ on $\td{\cal A}$ by
$$[u_{[0]},v_{[0]}]=[\ptl_1(u)\ptl_2(v)-\ptl_1(v)\ptl_2(u)+\ves(u\ptl_1(v)-\ptl_1(u)v)]_{[0]},\eqno(1.24)$$
$$[u_{[1]},v_{[1]}]=(x^{\al_0,\vec{0}}uv)_{[0]},\eqno(1.25)$$
\begin{eqnarray*}& &[u_{[0]},v_{[1]}]=-[v_{[1]},u_{[0]}]=(\ptl_1(u)\ptl_2(v)-\ptl_1(v)\ptl_2(u))_{[1]}\\& &+\ves\left(u\ptl_1(v)-{1\over 2}\ptl_1(u)v\right)_{[1]}+{1\over 2}(\vf_2(\al_0)\ptl_1(u)v+\vf_1(\al_0)(\ves uv-\ptl_2(u)v)_{[1]}\hspace{1.6cm}(1.26)\end{eqnarray*}
for $u,v\in{\cal A}_2$. We have verified that $(\td{\cal A},[\cdot,\cdot])$ forms a Lie superalgebra. In fact, the Lie superalgebra $({\cal L},[\cdot,\cdot])$ in Theorem 5.3 of [X2] is a special case of $(\td{\cal A},[\cdot,\cdot])$. We shall not present the tedious lengthy verification of the Jacobi identity for $(\td{\cal A},[\cdot,\cdot])$ in this paper. Up to equivalence of constructions, we only need to consider the cases when $\ves=0,1$. The special case when $\ves=0$ and $\al_0=0$ was included in the supersymmetric extension of the Lie algebras of Cartan type H in Section 6.3 of [X6]. Set
$$\td{\cal B}=\td{\cal A}_0+\td{\cal B}_1,\qquad \td{\cal B}_1=[\td{\cal A}_0,\td{\cal A}_1].\eqno(1.27)$$
Then $\td{\cal B}$ forms a subalgebra. Moreover, we set
$$\td{\cal C}=\left\{\begin{array}{ll}\td{\cal B}/\Bbb{F}1_{[0]}&\mbox{if}\;\;\ves=0,\\\td{\cal B}/\Bbb{F}(x^{\sgm_1,\vec{0}})_{[0]}&\mbox{if}\;\;\ves=1,\end{array}\right.\eqno(1.28) $$
where $\sgm_1\in\kn_{\vf_1}$ and $\vf_2(\sgm_1)=1$ as in Class I. Now we allow the case $\ves=1, \;\vf_2\equiv 0$ and ${\cal J}_2=\{0\}$; equivalently, we only assume that (1.6) only holds for $p=1$ when $\ves=1$.
 
\psp

{\bf Theorem 3}. {\it The pair} $(\td{\cal C},[\cdot,\cdot])$ {\it forms a simple Lie superalgebra  if} $\ves=1$ {\it or}  $\ves=0,\;\al_0\neq 0$, and
$$\kn_{\vf_1}\not\subset\kn_{\vf_2}\;\;\mbox{if}\;\;\vf_2\not\equiv 0\;\;\mbox{and}\;\;\kn_{\vf_2}\not\subset\kn_{\vf_1}\;\;\mbox{if}\;\;\vf_1\not\equiv 0.\eqno(1.29)$$
{\it  Moreover, the codimension of} $\td{\cal B}$ {\it in} $\td{\cal A}$ {\it is at most 1 under the above condition (1.29).} 
\psp

The above class of simple Lie superalgebras include the centerless super Virasoro algebra as a special case.

The following fact in linear algebra will be used in the proofs of Theorem 1 and Theorem 2. 
\psp

{\bf Lemma 4}. {\it Let} $T$ {\it be a linear transformation on a vector space} $U$ {\it and let} $U_1$ {\it be a subspace of} $U$ {\it such that} $T(U_1)\subset U_1$. {\it Suppose that} $u_1,u_2,...,u_n$ {\it are eigenvectors of} $T$ {\it corresponding to different eigenvalues. If}  $\sum_{j=1}^nu_j\in U_1$, {\it then} $u_1,u_2,...,u_n\in U_1$.
\pse

We shall present the proofs and related special examples of simple Lie algebras of Theorems 1, 2 and 3 in Sections 2, 3 and 4, respectively.

\section{Proof of Theorem 1 and Examples}

In this section, we shall first give the proof of Theorem 1 and then present some related examples of simple Lie algebras.
 \psp

\subsection{Proof of Theorem 1}

First we have
\begin{eqnarray*}[x^{\al,\vec{i}},x^{\be,\vec{j}}]_1&=&(\vf_1(\al)\vf_2(\be)-\vf_1(\be)\vf_2(\al)+\vf_1(\be-\al))x^{\al+\be,\vec{i}+\vec{j}}\\& &+(j_2\vf_1(\al)-i_2\vf_1(\be))x^{\al+\be,\vec{i}+\vec{j}-1_{[2]}}+(i_1j_2-i_2j_1)x^{\al+\be,\vec{i}+\vec{j}-1_{[1]}-1_{[2]}}\\& &+(i_1\vf_2(\be)-j_1\vf_2(\al)+j_1-i_1)x^{\al+\be,\vec{i}+\vec{j}-1_{[1]}}\hspace{4.6cm}(2.1)\end{eqnarray*}
for $(\al,\vec{i}),(\be,\vec{j})\in \G\times\vec{\cal J}$ by (1.13). In particular,
$$[1,x^{\be,\vec{j}}]_1=\vf_1(\be)x^{\be,\vec{j}}+j_1x^{\be,\vec{j}-1_{[1]}}\qquad\mbox{for}\;\;(\be,\vec{j})\in\G\times\vec{\cal J}.\eqno(2.2)$$
 Set
$$\G'=\{\vf_1(\al)\mid \al\in \G\}.\eqno(2.3)$$
For $\lmd\in\G'$, we define
$${\cal A}_{(\lmd)}=\{u\in{\cal A}_2\mid (\ad_1-\lmd)^m(u)=0\;\mbox{for some}\;m\in\Bbb{N}\}.\eqno(2.4)$$
By (2.2), we have
$${\cal A}_2=\bigoplus_{\lmd\in\G'}{\cal A}_{(\lmd)}.\eqno(2.5)$$

Recall that $\sgm_i\in\kn_{\vf_1}$ satisfying $\vf_2(\sgm_i)=i$ for $i=1,2$. 
\psp

{\it Proof of the First Statement in Theorem 1}
\psp

Let $I$ be any ideal of ${\cal A}_2$ that strictly contains $\Bbb{F}x^{\sgm_1,\vec{0}}$. To prove the first statement in Theorem 1 is equivalent to proving $I={\cal A}_2$. 

By (2.2) and (2.5), we have 
$$I=\bigoplus_{\lmd\in\G'}I_{\lmd},\qquad I_{\lmd}=I\bigcap {\cal A}_{(\lmd)}.\eqno(2.6)$$
For any $(\be,\vec{j})\in \G\times\vec{\cal J}$,
$$[x^{-\be,\vec{0}},x^{\be,\vec{j}}]_1=2\vf_1(\be)x^{0,\vec{j}}+j_1(\vf_2(\be)+1)x^{0,\vec{j}-1_{[1]}}-j_2\vf_1(\be)x^{0,\vec{j}-1_{[2]}}\eqno(2.7)$$
by (2.1) (cf. (1.9)). If $I_{\lmd}\neq\{0\}$ for some $0\neq \lmd\in\G'$, we pick any $0\neq u=\sum_{(\al,\vec{i})\in \G\times\vec{\cal J}}c_{\al,\vec{i}}x^{\al,\vec{i}}\\ \in I_{\lmd}$, where $c_{\al,\lmd}\in\Bbb{F}$. Assume $c_{\be,\vec{j}}\neq 0$ for some $(\be,\vec{j})\in \G\times\vec{\cal J}$. Fixing $\be$ and considering the nonzero terms $c_{\be,\vec{l}}x^{\be,\vec{l}}$ in $u$ with the largest value $l_1+l_2$, we have
$$[x^{-\be,\vec{0}},u]_1\in I_0\setminus \Bbb{F}x^{\sgm_1,\vec{0}}\eqno(2.8)$$
by (2.7) and the fact that $\sgm_1\neq 0$. Thus (2.6) implies that  $I_0$ strictly contains $\Bbb{F}x^{\sgm_1,\vec{0}}$. 

Equation (2.2) imply 
$$[1,x^{\be,\vec{j}}]_1=j_1x^{\be,\vec{j}-1_{[1]}}\qquad\mbox{for}\;\;(\be,\vec{j})\in\kn_{\vf_1}\times\vec{\cal J}.\eqno(2.9)$$
Moreover, we have
$$[x^{-\sgm_1,\vec{0}},x^{\sgm_1,1_{[1]}}]_1=2x^{0,\vec{0}}\eqno(2.10)$$
by (2.7). Set
$$I_{[0]}=I\bigcap(\sum_{(\al,j)\in\kn_{\vf_1}\times{\cal J}_2}\Bbb{F}x^{\al,j_{[2]}}).\eqno(2.11)$$
By (2.9), (2.10) and repeated acting $\ad_1$ on $I_0$, we can prove that 
$$I_{[0]}\;\;\mbox{strictly contains}\;\; \Bbb{F}x^{\sgm_1,\vec{0}}.\eqno(2.12)$$

For $k\in\Bbb{N}$, we let
$${\cal A}_{\la k\ra}=\mbox{span}\{x^{\al,\vec{i}}\mid (\al,\vec{i})\in \G\times \vec{\cal J},\;i_2\leq k\}.\eqno(2.13)$$
Then
$${\cal A}_2=\bigcup_{k=0}^{\infty}{\cal A}_{\la k\ra}.\eqno(2.14)$$
Set
$$\hat{k}=\min\{k\in\Bbb{N}\mid I_{[0]}\bigcap {\cal A}_{\la k\ra}\setminus \Bbb{F}x^{\sgm_1,\vec{0}}\neq\emptyset\}.\eqno(2.15)$$
For any $u\in (I_{[0]}\bigcap {\cal A}_{\la \hat{k}\ra})\setminus \Bbb{F}x^{\sgm_1,\vec{0}}$, we write it as
$$u=\sum_{\al\in\kn_{\vf_1},(\al,\hat{k})\neq (\sgm_1,\vec{0})}a_{\al}x^{\al,\hat{k}_{[2]}}+u',\qquad a_{\al}\in\Bbb{F},\;u'\in{\cal A}_{\la \hat{k}-1\ra}+\Bbb{F}x^{\sgm_1,\vec{0}}\eqno(2.16)$$
and define
$$\imath(u)=|\{\al\mid a_{\al}\neq 0\}|.\eqno(2.17)$$
By (2.12), $\imath(u)>0$. Set
$$\imath=\min\{\imath(v)\mid 0\neq v\in (I_{[0]}\bigcap {\cal A}_{\la \hat{k}\ra})\setminus \Bbb{F}x^{\sgm_1,\vec{0}}\}.\eqno(2.18)$$
Choose $0\neq u\in (I_{[0]}\bigcap {\cal A}_{\la \hat{k}\ra})\setminus \Bbb{F}x^{\sgm_1,\vec{0}}$ such that $\imath(u)=\imath$. Write $u$ as in (2.16).

If ${\cal J}_1=\Bbb{N}$, then
$$I_{[0]}\bigcap {\cal A}_{\la \hat{k}\ra}\ni [x^{0,1_{[1]}},u]_1\equiv \sum_{\al\in\kn_{\vf_1},(\al,\hat{k})\neq (\sgm_1,\vec{0})}a_{\al}(\vf_2(\al)-1)x^{\al,\hat{k}_{[2]}}\;\;(\mbox{mod}\;{\cal A}_{\la \hat{k}-1\ra})\eqno(2.19)$$
by (2.1). By the minimality of $\imath(u)$ (cf. (2.18)) and Lemma 4,
$$\vf_2(\al)=\vf_2(\be)\qquad\mbox{whenever}\;\;a_{\al}a_{\be}\neq 0.\eqno(2.20)$$
Assume that ${\cal J}_1=\{0\}$. By (1.6), $\vf_1\not\equiv 0$. For any $\tau\in\G\setminus\kn_{\vf_1}$, we have
\begin{eqnarray*}& &I_{[0]}\bigcap {\cal A}_{\la \hat{k}\ra}\ni[[x^{\tau,\vec{0}},u]_1,x^{-\tau,\vec{0}}]_1\equiv\\& & \sum_{\al\in\kn_{\vf_1},(\al,\hat{k})\neq (\sgm_1,\vec{0})}a_{\al}\vf_1(\tau)^2(\vf_2(\al)-1)(\vf_2(\al)-2)x^{\al,\hat{k}_{[2]}}\;\;(\mbox{mod}\;{\cal A}_{\la \hat{k}-1\ra}).\hspace{1.6cm}(2.21)\end{eqnarray*}
Again by the minimality of $\imath(u)$ (cf. (2.18)) and Lemma 4,
$$(\vf_2(\al)-1)(\vf_2(\al)-2)=(\vf_2(\be)-1)(\vf_2(\be)-2)\qquad\mbox{whenever}\;\;a_{\al}a_{\be}\neq 0,\eqno(2.22)$$
which can be written as
$$(\vf_2(\al)-\vf_2(\be))(\vf_2(\al)+\vf_2(\be)-3)=0\qquad\mbox{whenever}\;\;a_{\al}a_{\be}\neq 0.\eqno(2.23)$$
Thus we obtain
$$\vf_2(\al)=\vf_2(\be)\;\;\mbox{or}\;\;\vf_2(\al)+\vf_2(\be)=3\qquad\mbox{whenever}\;\;a_{\al}a_{\be}\neq 0.\eqno(2.24)$$
Assume that there exists $a_{\al}a_{\be}\neq 0$ such that $\vf_2(\al)\neq\vf_2(\be)$ and $\vf_2(\al)+\vf_2(\be)=3$. Without loss of generality, we can assume $\vf_2(\be)\neq 0, 1$ because $\vf_2(\al)+\vf_2(\be)=3$. For any $\tau\in\G\setminus\kn_{\vf_1}$, we have
\begin{eqnarray*}&  &I_{[0]}\bigcap {\cal A}_{\la \hat{k}\ra}\ni [[x^{\al+\tau,\vec{0}},u]_1,x^{-\tau,\vec{0}}]_1\equiv \sum_{\gm\in\kn_{\vf_1},(\gm,\hat{k})\neq (\sgm_1,\vec{0})}\\& &a_{\gm}\vf_1(\tau)^2(\vf_2(\gm)-1)(\vf_2(\al+\gm)-2)x^{\al+\gm,\hat{k}_{[2]}}\;\;(\mbox{mod}\;{\cal A}_{\la \hat{k}-1\ra}).\hspace{3.8cm}(2.25)\end{eqnarray*}
Since $\vf_2$ is a group homomorphism, $\vf_2(\al+\be)=\vf_2(\al)+\vf_2(\be)=3$. Hence
$$a_{\be}\vf_1(\tau)^2(\vf_2(\be)-1)(\vf_2(\al+\be)-2)x^{\al+\be,\hat{k}_{[2]}}\neq 0.\eqno(2.26)$$
Thus
$$0<\imath([[x^{\al+\tau,\vec{0}},u]_1,x^{-\tau,\vec{0}}]_1)\leq \imath(u)=\imath\eqno(2.27)$$
because $\vf_2(\al+\be)=3\neq 1=\vf_2(\sgm_1)$.
By the minimality of $\imath(u)$ (cf. (2.18)), we have
$$\imath([[x^{\al+\tau,\vec{0}},u]_1,x^{-\tau,\vec{0}}]_1)=\imath(u),\eqno(2.28)$$
which implies
$$a_{\al}\vf_1(\tau)^2(\vf_2(\al)-1)(\vf_2(2\al)-2)\neq 0\eqno(2.29)$$
by (2.25).
So $\vf_2(\al)\neq 1$. Similarly, we have
\begin{eqnarray*}& &I_{[0]}\bigcap {\cal A}_{\la \hat{k}\ra}\ni [[x^{\be+\tau,\vec{0}},u]_1,x^{-\tau,\vec{0}}]_1\equiv\sum_{\gm\in\kn_{\vf_1},(\gm,\hat{k})\neq (\sgm_1,\vec{0})}\\& & a_{\gm}\vf_1(\tau)^2(\vf_2(\gm)-1)(\vf_2(\be+\gm)-2)x^{\be+\gm,\hat{k}_{[2]}}\;\;(\mbox{mod}\;{\cal A}_{\la \hat{k}-1\ra}),\hspace{3.6cm}(2.30)\end{eqnarray*}
$$\imath([[x^{\be+\tau,\vec{0}},u]_1,x^{-\tau,\vec{0}}]_1)=\imath(u)\eqno(2.31)$$
and
$$a_{\al}\vf_1(\tau)^2(\vf_2(\al)-1)(\vf_2(\al+\be)-2)a_{\be}\vf_1(\tau)^2(\vf_2(\be)-1)(\vf_2(2\be)-2)\neq 0.\eqno(2.32)$$
But
$$\vf_2(\al+\be)\neq \vf_2(2\be),\;\;\vf_2(\al+\be)+\vf_2(2\be)=2\vf_2(\be)+3\neq 3,\eqno(2.33)$$
which contradicts  (2.24) if we replace $u$ by $[[x^{\be+\tau,\vec{0}},u]_1,x^{-\tau,\vec{0}}]_1$. Thus we always have (2.20). Therefore, by (1.7), we have
$$\al=\be\qquad\mbox{whenever}\;\;a_{\al}a_{\be}\neq 0.\eqno(2.34)$$
Hence we can assume
$$u=x^{\al,\hat{k}_{[2]}}+u',\;\; \al\in\kn_{\vf_1},\;u'\in{\cal A}_{\la \hat{k}-1\ra}+\Bbb{F}x^{\sgm_1,\vec{0}}\eqno(2.35)$$
and $\al\neq\sgm_1$ if $\hat{k}=0$.

Assume that $\hat{k}>0$ and  $\vf_2(\al)=1$. If ${\cal J}_1=\Bbb{N}$, then
$$I_{[0]}\bigcap{\cal A}_{\la \hat{k}-1\ra}\ni [x^{-\al,1_{[1]}},u]=\hat{k}x^{0,(\hat{k}-1)_{[2]}}+\sum_{0\neq\be\in\kn_{\vf_1}}c_{\be}x^{\be,(\hat{k}-1)_{[2]}}\;\;(\mbox{mod}\;{\cal A}_{\la \hat{k}-2\ra})\eqno(2.36)$$
with $c_{\be}\in\Bbb{F}$, which implies $(I_{[0]}\bigcap{\cal A}_{\la \hat{k}-1\ra})\setminus \Bbb{F}x^{\sgm_1,\vec{0}}\neq\emptyset$. We get a contradiction to (2.15). If  ${\cal J}_1=\{0\}$, then $\vf_1\not\equiv 0$ by (1.6) and for any $\tau\in\G\setminus\kn_{\vf_1}$,
\begin{eqnarray*}&&I_{[0]}\bigcap{\cal A}_{\la \hat{k}-1\ra}\ni[[u,x^{-\al-\tau,\vec{0}}]_1,x^{\tau,\vec{0}}]_1\\&\equiv& 2\hat{k}\vf_1(\tau)^2x^{0,(\hat{k}-1)_{[2]}}+\sum_{0\neq\be\in\kn_{\vf_1}}c_{\be}x^{\be,(\hat{k}-1)_{[2]}}\;\;(\mbox{mod}\;{\cal A}_{\la \hat{k}-2\ra}),\hspace{3.6cm}(2.37)\end{eqnarray*}
where $c_{\be}\in\Bbb{F}$. Again we get a contradiction to (2.15). 

Suppose that $\hat{k}>0$ and  $\vf_2(\al)\neq 1$. If $\vf_1\not\equiv 0$, then for any $\tau\in\G\setminus\kn_{\vf_1}$, we have
\begin{eqnarray*}&&I_{[0]}\bigcap{\cal A}_{\la \hat{k}-1\ra}\ni[x^{-2\tau,\vec{0}},[x^{\tau,\vec{0}},[u,x^{-\al+\tau,\vec{0}}]_1]_1]_1\\&\equiv& 4\hat{k}\vf_1(\tau)^3(1-\vf_2(\al))x^{0,(\hat{k}-1)_{[2]}}+\sum_{0\neq\be\in\kn_{\vf_1}}c_{\be}x^{\be,(\hat{k}-1)_{[2]}}\;\;(\mbox{mod}\;{\cal A}_{\la \hat{k}-2\ra})\hspace{1.9cm}(2.38)\end{eqnarray*}
with $c_{\be}\in\Bbb{F}$, which contradicts  (2.15). If $\vf_1\equiv 0$, then ${\cal J}_1=\Bbb{N}$ by (1.6). Moreover,
$$[x^{-\al,1_{[1]}},u]\equiv (\vf_2(\al)-1)x^{0,\hat{k}_{[2]}}\;\;(\mbox{mod}\;{\cal A}_{\la \hat{k}-1\ra})\eqno(2.39)$$
by (2.1) and the facts that $\vf_1\equiv 0$. Replacing $u$ by $(\vf_2(\al)-1)^{-1}[x^{-\al,1_{[1]}},u]$,  we can assume that $\al=0$. Furthemore,
$$[x^{0,1_{[1]}},u]+u=\ptl_2(u)\in I_{[0]}\bigcap{\cal A}_{\la \hat{k}-1\ra}\setminus\Bbb{F}x^{\sgm_1,\vec{0}}\eqno(2.40)$$
by (1.12) and (1.13), which contradicts (2.15).

Thus we have $\hat{k}=0$. So $x^{\al,\vec{0}}\in I$ for some $\al\in\kn_{\vf_1}$ such that $\vf_2(\al)\neq 1$. Assume that $\al\neq 0$. So $\vf_2(\al)\neq 0$ by (1.7). If ${\cal J}_1=\Bbb{N}$, then
$$[x^{-\al,1_{[1]}},x^{\al,\vec{0}}]_1=(\vf_2(\al)-1)x^{0,\vec{0}}\in I.\eqno(2.41)$$
If ${\cal J}_1=\{0\}$, then $\vf_1\not\equiv 0$ by (1.6).  Picking any $\tau\in\G\setminus\kn_{\vf_1}$, we have
$$[[x^{\al,\vec{0}},x^{-\al-\tau,\vec{0}}]_1,x^{\tau,\vec{0}}]_1=2\vf_1(\tau)^2(\vf_2(\al)-1)x^{0,\vec{0}}\in I.\eqno(2.42)$$
Thus we always have $1=1_{\cal A}=x^{0,\vec{0}}\in I$.
 Moreover, we can prove $I={\cal A}_2$ if ${\cal J}_1=\Bbb{N}$ by (2.2) and induction $j_1$. Assume that ${\cal J}_1=\{0\}$. So $\vf_1\not\equiv 0$ by (1.6). Then (2.2) implies
$$x^{\al,j_{[2]}}\in I\qquad\mbox{for}\;\;(\al,j)\in\G\times{\cal J}_2,\;\vf_1(\al)\neq 0.\eqno(2.43)$$
For any $(\be,j)\in \kn_{\vf_1}\times{\cal J}_2$ and $\tau\in\G\setminus\kn_{\vf_1}$, we have:
$$[x^{\tau,\vec{0}},x^{-\tau+\be,j_{[2]}}]_1=\vf_1(\tau)((\vf_2(\be)-2)x^{\be,j_{[2]}}+jx^{\be,(j-1)_{[2]}})\in I.\eqno(2.44)$$
By (2.44) and  induction $j$, we can prove $I={\cal A}$ if ${\cal J}_2=\Bbb{N}$. If ${\cal J}_2=\{0\}$, then (2.44) shows
$$ x^{\be,\vec{0}}\in I\qquad \mbox{for}\;\;\sgm_2\neq \be\in\kn_{\vf_1}.\eqno(2.45)$$
Thus $I={\cal A}$ by the assumption of the first statement of Theorem 1.
\psp

{\it Proof of the Second Statement in Theorem 1}
\psp

Set 
$$\hat{\cal B}_2=\mbox{span}\:\{x^{\al,\vec{0}}\mid\sgm_2\neq\al\in\G\}.\eqno(2.46)$$
Then (2.2) and (2.44) shows
$$\hat{\cal B}_2=[{\cal A}_2,{\cal A}_2]_1.\eqno(2.47)$$
Replacing ${\cal A}$ by $\hat{\cal B}_2$ and $\G$ by $\G\setminus\{\sgm_2\}$ in the proof of the first statement, we obtain (2.45), which implies the simplicity of
$${\cal B}_2^{(1)}=[{\cal B}_2,{\cal B}_2]=\hat{\cal B}_2/\Bbb{F}x^{\sgm_1,\vec{0}}.\eqno(2.48)$$

This completes the proof of Theorem 1.

\subsection{Related Examples of Simple Lie Algebras}

Let $\Bbb{R}$ be the field of real numbers. By Theorem 1, we obtain the following examples of simple Lie algebras $({\cal B},[\cdot,\cdot])$, where $f_t$ denotes the partial derivative of a polynomial $f$ with respect to the variable $t$.
\psp

{\bf Example 2.1}. The space ${\cal B}=\Bbb{R}[t_1,t_2]$ with the Lie bracket:
$$[f,g]=f_{t_1}(g_{t_2}-g)+(f-f_{t_2})g_{t_1}\qquad\mbox{for}\;\;f,g\in \Bbb{R}[t_1,t_2].\eqno(2.49)$$

{\bf Example 2.2}. The space ${\cal B}=\Bbb{R}[t_1,t_2, t_1^{-1}]$ with the Lie bracket:
$$[f,g]=t_1[f_{t_1}(g_{t_2}-g)+(f-f_{t_2})g_{t_1}].\eqno(2.50)$$
The space ${\cal B}=\Bbb{R}[t_1,t_2, t_1^{-1}]/\Bbb{R}t_1$  and its Lie bracket is induced by:
$$[f,g]=f_{t_2}(t_1g_{t_1}-g)+(f-t_1f_{t_1})g_{t_2}\eqno(2.51)$$
for $f,g\in \Bbb{R}[t_1,t_2,t_1^{-1}]$.
\pse

{\bf Example 2.3}. The space ${\cal B}=\Bbb{R}[t_1,t_2,t_3, t_1^{-1}]$ with the Lie bracket:
$$[f,g]=(t_1f_{t_1}+f_{t_3})(g_{t_2}-g)+(f-f_{t_2})(t_1g_{t_1}+g_{t_3}).\eqno(2.52)$$
The space ${\cal B}=\Bbb{R}[t_1,t_2,t_3, t_1^{-1}]/\Bbb{R}t_1$ with the Lie bracket induced by:
$$[f,g]=f_{t_2}(t_1g_{t_1}+g_{t_3}-g)+(f-t_1f_{t_1}-f_{t_3})g_{t_2}\eqno(2.53)$$
for $f,g\in \Bbb{R}[t_1,t_2,t_3,t_1^{-1}]$.
\pse

{\bf Example 2.4}. The space ${\cal B}=\Bbb{R}[t_1,t_2,t_3, t_1^{-1},t_2^{-1}]/\Bbb{R}t_2$ with the Lie bracket induced by:
$$[f,g]=(t_1f_{t_1}+f_{t_3})(t_2g_{t_2}-g)+(f-t_2f_{t_2})(t_1g_{t_1}+g_{t_3})\eqno(2.54)$$
or
$$[f,g]=t_1[f_{t_1}(t_2g_{t_2}+g_{t_3}-g)+(f-t_2f_{t_2}-f_{t_3})g_{t_1}]\eqno(2.55)$$
for $f,g\in \Bbb{R}[t_1,t_2,t_3,t_1^{-1},t_2^{-1}]$.
\pse

{\bf Example 2.5}. Let $m$ be a positive integer. Let ${\cal A}$ be the subalgebra of 
$$\Bbb{R}[t_1^{1/m},t_2^{1/m},t_3,t_4, t_1^{-1/m},t_2^{-1/m}]\eqno(2.56)$$
 generated by 
$$\{t_1,t_2,t_3,t_4,(t_1t_2)^{1/m},t_1^{-1},t_2^{-1}\}.\eqno(2.57)$$
In particular, ${\cal A}=\Bbb{R}[t_1,t_2,t_3,t_4, t_1^{-1},t_2^{-1}]$ when $m=1$. The space ${\cal B}={\cal A}/\Bbb{R}t_2$ and its Lie bracket is induced by 
$$[f,g]=(t_1f_{t_1}+f_{t_3})(t_2g_{t_2}-g_{t_4}-g)+(f-t_2f_{t_2}-f_{t_4})(t_1g_{t_1}+g_{t_3})\eqno(2.58)$$
for $f,g\in{\cal A}$.

\section{Proof of Theorem 2 and Examples}

In this section, we shall first give the proof of Theorem 2 and then present some related examples of simple Lie algebras.

\subsection{Proof of Theorem 2}

 First we have
\begin{eqnarray*}& &[x^{\al,\vec{i}},x^{\be,\vec{j}}]_2=\\& &(\vf_1(\al)\vf_2(\be)-\vf_1(\be)\vf_2(\al))x^{\al+\be+\al_0,\vec{i}+\vec{j}}+(i_1\vf_2(\be)-j_1\vf_2(\al))x^{\al+\be+\al_0,\vec{i}+\vec{j}-1_{[1]}}\\& &+(j_2\vf_1(\al)-i_2\vf_1(\be)x^{\al+\be+\al_0,\vec{i}+\vec{j}-1_{[2]}}+(i_1j_2-i_2j_1))x^{\al+\be+\al_0,\vec{i}+\vec{j}-1_{[1]}-1_{[2]}}\\& &+(\vf_3(\al+\al_0)\vf_4(\be+\al_0)-\vf_3(\be+\al_0)\vf_4(\al+\al_0))x^{\al+\be,\vec{i}+\vec{j}}\\& &+(i_3\vf_4(\be+\al_0)-j_3\vf_4(\al+\al_0))x^{\al+\be,\vec{i}+\vec{j}-1_{[3]}}+(i_3j_4-i_4j_3)x^{\al+\be,\vec{i}+\vec{j}-1_{[3]}-1_{[4]}}\\& &+(j_4\vf_3(\al+\al_0)-i_4\vf_3(\be+\al_0))x^{\al+\be,\vec{i}+\vec{j}-1_{[4]}},\hspace{6.5cm}(3.1)\end{eqnarray*}
for $(\al,\vec{i}),(\be,\vec{j})\in \G\times\vec{\cal J}$ by (1.18).
Recall the element $\sgm\in\kn_{\vf_1}\bigcap\kn_{\vf_2}$ satisfying $\vf_3(\sgm+\al_0)=\vf_4(\sgm+\al_0)=0$. Moreover, we treat $x^{\sgm,\vec{0}}$ as $0\in {\cal A}_4$ when such $\sgm$ does not exist. 
\psp

{\it Proof of the First Statement in Theorem 2}.
\psp

 Let $I$ be any ideal of ${\cal A}_4$ that strictly contains $\Bbb{F}x^{\sgm,\vec{0}}$. To prove the first statement in Thorem 2 is equivalent to proving $I={\cal A}_4$.
\psp

{\it Step 1}. $1\in I$.
\psp

For any $\vec{i}\in\vec{\cal J}$, we define
$$|\vec{i}|=\sum_{p=1}^4i_p.\eqno(3.2)$$
Set 
$${\cal A}_{4,k}=\mbox{span}\:\{x^{\al,\vec{i}}\mid \al\in\G,\;\vec{i}\in\vec{\cal J},\;|\vec{i}|\leq k\}\qquad\mbox{for}\;\;k\in\Bbb{N}.\eqno(3.3)$$
For convenience, we let
$${\cal A}_{4,-1}=\emptyset.\eqno(3.4)$$
Moreover, we define 
$$\hat{k}=\mbox{min}\;\{k\in\Bbb{N}\mid ({\cal A}_{4,k}\bigcap I)\setminus \Bbb{F}x^{\sgm,\vec{0}}\neq \emptyset\}.\eqno(3.5)$$
For any $u\in ({\cal A}_{4,\hat{k}}\bigcap I)\setminus \Bbb{F}x^{\sgm,\vec{0}}$, we write:
$$u=u_{ld}+\td{u}\qquad\mbox{with}\;\;\td{u}\in {\cal A}_{4,\hat{k}-1}+\Bbb{F}x^{\sgm,\vec{0}}\eqno(3.6)$$
and
$$u_{ld}=\sum_{(\sgm,\vec{0})\neq (\al,\vec{i})\in\G\times\vec{\cal J},\:|\vec{i}|=\hat{k}}a_{\al,\vec{i}}x^{\al,\vec{i}},\qquad\;\;a_{\al,\vec{i}}\in\Bbb{F}\eqno(3.7)$$
and  define
$$\flat(u)=|\{\al\in\G\mid a_{\al,\vec{i}}\neq 0\;\mbox{for}\; \vec{i}\in\vec{\cal J},\;|\vec{i}|=\hat{k}\}|.\eqno(3.8)$$
Furthermore, we set
$$\flat=\min\{\flat(v)\mid v\in ({\cal A}_{4,\hat{k}}\bigcap I)\setminus\Bbb{F}x^{\sgm,\vec{0}}\}.\eqno(3.9)$$

Let $ u\in ({\cal A}_{4,\hat{k}}\bigcap I)\setminus\Bbb{F}x^{\sgm,\vec{0}}$ such that $\flat(u)=\flat$. Write $u$ as in (3.6) and (3.7). We set
$$1'=2,\;\;2'=1,\;\;3'=4,\;\;4'=3,\;\;\es(1)=\es(3)=1,\;\; \es(2)=\es(4)=-1.\eqno(3.10)$$

For $p\in \{1,2\}$ and any $\tau\in\bigcap_{p\neq q\in\ol{1,4}}\kn_{\vf_q}$,
we have
\begin{eqnarray*}& &[x^{-\tau-\al_0,\vec{0}},[u,x^{\tau-\al_0,\vec{0}}]_2]_2\\&\equiv& \sum_{(\sgm,\vec{0})\neq (\al,\vec{i})\in\G\times\vec{\cal J},\:|\vec{i}|=\hat{k}}[(\vf_p(\tau)-\vf_p(\al_0))\vf_{p'}(\al)+\vf_{p'}(\al_0)\vf_p(\al)]\\& &
[(\vf_p(\tau)+\vf_p(\al_0))\vf_{p'}(\al)-\vf_{p'}(\al_0)(\vf_p(\al)+\vf_p(\tau)))]a_{\al,\vec{i}}x^{\al,\vec{i}}\;\;(\mbox{mod}\:{\cal A}_{\hat{k}-1})\hspace{1.2cm}(3.11)\end{eqnarray*}
by (3.1). Since $\vf_p(\tau)$ takes an infinite number of elements in $\Bbb{F}$ for  $\tau\in\bigcap_{p\neq q\in\ol{1,4}}\kn_{\vf_q}$ by (1.15), the coefficients of $\vf_p(\tau)^2$ show
$$\vf_{p'}(\al)\vf_{p'}(\al-\al_0)=\vf_{p'}(\be)\vf_{p'}(\be-\al_0)\;\;\;\mbox{whenever}\;\;a_{\al,\vec{i}}a_{\be,\vec{j}}\neq 0\eqno(3.12)$$
by the minimality of $\flat(u)$ (cf. (3.9)) and Lemma 4. Moreover, (3.12) is equivalent to
$$\vf_{p'}(\al)=\vf_{p'}(\be)\;\;\mbox{or}\;\;\vf_{p'}(\al+\be-\al_0)=0\;\;\;\mbox{whenever}\;\;a_{\al,\vec{i}}a_{\be,\vec{j}}\neq 0.\eqno(3.13)$$
Assume that there exist $\al,\be\in\G$ such that $\vf_{p'}(\al)\neq \vf_{p'}(\be),\;\vf_{p'}(\al+\be-\al_0)=0$ and $a_{\al,\vec{i}}a_{\be,\vec{j}}\neq 0$. We may assume that $\vf_{p'}(\al)\neq 0$. Choose
$$\tau\in\bigcap_{p\neq q\in\ol{1,4}}\kn_{\vf_q}\setminus\kn_{\vf_p},\;\;\tau'\in\bigcap_{p'\neq q\in\ol{1,4}}\kn_{\vf_q}\setminus\kn_{\vf_{p'}}\eqno(3.14)$$
such that
$$\vf_{p'}(\al+\tau'+\al_0)\neq 0,\;\;\vf_p(\tau-\al_0)\vf_{p'}(\al)\neq \vf_p(\al)\vf_{p'}(\tau'-\al_0)\eqno(3.15)$$
by (1.15). We have
\begin{eqnarray*}[x^{\tau+\tau'-\al_0,\vec{0}},u]_2&\equiv& \sum_{(\sgm,\vec{0})\neq (\gm,\vec{l})\in\G\times\vec{\cal J},\:|\vec{l}|=\hat{k}}\es_pa_{\gm,\vec{l}}[\vf_p(\tau-\al_0)\vf_{p'}(\gm)\\& &-\vf_p(\gm)\vf_{p'}(\tau'-\al_0)]x^{\gm+\tau+\tau',\vec{l}}\;\;(\mbox{mod}\;{\cal A}_{\hat{k}-1})\hspace{4cm}(3.16)\end{eqnarray*}
by (3.1). Since $\flat(u)$ is minimum, $\al+\tau+\tau'\neq \sgm$ due to $\vf_{p'}(\al+\tau'+\al_0)\neq 0$, and 
$$\es_p a_{\al,\vec{i}}(\vf_p(\tau-\al_0)\vf_{p'}(\al)-\vf_p(\al)\vf_{p'}(\tau'-\al_0))\neq 0,\eqno(3.17)$$
we have
$$\es_pa_{\be,\vec{j}}(\vf_p(\tau-\al_0)\vf_{p'}(\be)-\vf_p(\be)\vf_{p'}(\tau'-\al_0))\neq 0.\eqno(3.18)$$
But 
$$\vf_{p'}(\al+\tau+\tau')\neq \vf_{p'}(\be+\tau+\tau'),\;\;\vf_{p'}((\al+\tau+\tau')+(\be+\tau+\tau')-\al_0)=2\vf_{p'}(\tau')\neq 0,\eqno(3.19)$$
which contradicts to (3.13) with $u$ replaced by $[x^{\tau'-\al_0,\vec{0}},u]_2$. Thus the first equation in (3.13) holds.

For $p\in \{3,4\}$ and any $\tau\in\bigcap_{p\neq q\in\ol{1,4}}\kn_{\vf_q}$,
we have
\begin{eqnarray*}& &[x^{-\tau,\vec{0}},[u,x^{\tau,\vec{0}}]_2]_2\\&\equiv& \sum_{(\sgm,\vec{0})\neq (\al,\vec{i})\in\G\times\vec{\cal J},\:|\vec{i}|=\hat{k}}[\vf_p(\tau)\vf_{p'}(\al+\al_0)+\vf_p(\al_0)\vf_{p'}(\al)-\vf_p(\al)\vf_{p'}(\al_0)]\\& &
[\vf_p(\tau)\vf_{p'}(\al+2\al_0)-\vf_p(\al_0)\vf_{p'}(\al)+\vf_p(\al)\vf_{p'}(\al_0)]
a_{\al,\vec{i}}x^{\al,\vec{i}}\;\;(\mbox{mod}\:{\cal A}_{\hat{k}-1})\hspace{1.2cm}(3.20)\end{eqnarray*}
by (3.1). Since $\vf_p(\tau)$ takes an infinite number of elements in $\Bbb{F}$ for  $\tau\in\bigcap_{p\neq q\in\ol{1,4}}\kn_{\vf_q}$ by (1.15), the coefficients of $\vf_p(\tau)^2$ show
$$\vf_{p'}(\al+\al_0)\vf_{p'}(\al+2\al_0)=\vf_{p'}(\be+\al_0)\vf_{p'}(\be+2\al_0)\;\;\;\mbox{whenever}\;\;a_{\al,\vec{i}}a_{\be,\vec{j}}\neq 0\eqno(3.21)$$
by the minimality of $\flat(u)$ (cf. (3.9)) and Lemma 4. Moreover, (3.21) is equivalent to
$$\vf_{p'}(\al)=\vf_{p'}(\be)\;\;\mbox{or}\;\;\vf_{p'}(\al+\be+3\al_0)=0\;\;\;\mbox{whenever}\;\;a_{\al,\vec{i}}a_{\be,\vec{j}}\neq 0.\eqno(3.22)$$
Assume that there exist $\al,\be\in\G$ such that $\vf_{p'}(\al)\neq \vf_{p'}(\be),\;\vf_{p'}(\al+\be+3\al_0)=0$ and $a_{\al,\vec{i}}a_{\be,\vec{j}}\neq 0$. We may assume that $\vf_{p'}(\al+\al_0)\neq 0$ because $\vf_{p'}(\al)\neq \vf_{p'}(\be)$. Choose $\tau$ and $\tau'$ as in (3.14) such that
$$\vf_{p'}(\al+\tau')\neq 0,\;\;\vf_p(\tau+\al_0)\vf_{p'}(\al+\al_0)\neq \vf_p(\al+\al_0)\vf_{p'}(\tau'+\al_0)\eqno(3.23)$$
by (1.15). We have
\begin{eqnarray*}[x^{\tau+\tau',\vec{0}},u]_2&\equiv& \sum_{(\sgm,\vec{0})\neq (\gm,\vec{l})\in\G\times\vec{\cal J},\:|\vec{l}|=\hat{k}}\es_{p}a_{\gm,\vec{l}}
[\vf_p(\tau+\al_0)\vf_{p'}(\gm+\al_0)\\& &-\vf_p(\gm+\al_0)\vf_{p'}(\tau'+\al_0)]x^{\gm+\tau+\tau',\vec{l}}\;\;(\mbox{mod}\;{\cal A}_{\hat{k}-1}).\hspace{3.5cm}(3.24)\end{eqnarray*}
by (3.1). As (3.17)-(3.19), we get a contradiction to (3.22) with $u$ replaced by $[x^{\tau+\tau',\vec{0}},u]_2$. Thus the first equation in (3.22) holds. By (1.7), we obtain
$$\al=\be\;\;\;\mbox{whenever}\;\;a_{\al,\vec{i}}a_{\be,\vec{j}}\neq 0.\eqno(3.25)$$

Let $a_{\al,\vec{i}}\neq 0$ be fixed.
\psp

{\it Case 1}. $\hat{k}=0$.
\psp

In this case, $x^{\al,\vec{0}}\in I$ for some $\sgm\neq\al\in\G$. Assume that $\al\neq 0$. Since
$$[x^{\al,\vec{0}},x^{-\al,\vec{0}}]_2=2(\vf_3(\al)\vf_4(\al_0)-\vf_3(\al_0)\vf_4(\al))x^{0,\vec{0}}\eqno(3.26)$$
by (3.1), then $1=1_{\cal A}=x^{0,\vec{0}}\in I$  if $\vf_3(\al)\vf_4(\al_0)\neq \vf_3(\al_0)\vf_4(\al)$. 

Assume that $\vf_3(\al)\vf_4(\al_0)=\vf_3(\al_0)\vf_4(\al)$ and $\vf_3(\al+\al_0)\neq 0$ or $\vf_4(\al+\al_0)\neq 0$. By (1.15) and (1.16), we can choose 
$$\tau\in\bigcap_{3\neq q\in\ol{1,4}}\kn_{\vf_q},\;\;\tau'\in\bigcap_{q=1}^3\kn_{\vf_q}\eqno(3.27)$$
such that
$$\vf_3(\tau)\vf_4(\al_0)\neq \vf_3(\al_0)\vf_4(\tau')\eqno(3.28)$$
and
$$\vf_3(\tau)\vf_4(\al+\al_0)\neq \vf_3(\al+\al_0)\vf_4(\tau').\eqno(3.29)$$
Then
$$[x^{\tau+\tau',\vec{0}},x^{\al,\vec{0}}]_2=(\vf_3(\tau)\vf_4(\al+\al_0)-\vf_3(\al+\al_0)\vf_4(\tau'))x^{\al+\tau+\tau',\vec{0}}\in I.\eqno(3.30)$$
Thus $x^{\al+\tau+\tau',\vec{0}}\in I$. Moreover,
\begin{eqnarray*}& &\vf_3(\al+\tau+\tau')\vf_4(\al_0)=(\vf_3(\al)+\vf_3(\tau))\vf_4(\al_0)\\& \neq &\vf_3(\al_0)(\vf_4(\al)+\vf_4(\tau'))=\vf_3(\al_0)\vf_4(\al+\tau+\tau').\hspace{5.4cm}(3.31)\end{eqnarray*}
So $1\in I$ by the arguments in the above paragraph. 

Suppose that  $\vf_3(\al+\al_0)=\vf_4(\al+\al_0)=0$. Since $\al\neq \sgm$, there exists $p\in\{1,2\}$ such that $\vf_p(\al)\neq 0$. Pick any $\tau'\in(\bigcap_{p'\neq q\in\ol{1,4}}\kn_{\vf_q})\setminus \kn_{\vf_{p'}}$. Then we have
$$[x^{\al,\vec{0}},x^{\tau',\vec{0}}]_2=\es_p\vf_p(\al)\vf_{p'}(\tau')x^{\al+\al_0+\tau',\vec{0}}\in I.\eqno(3.32)$$
By (1.16), there exists $q\in\{3,4\}$ such that $\vf_q(\al_0)\neq 0$. Thus
$$\vf_q((\al+\al_0+\tau')+\al_0)=\vf_q(\al_0)\neq 0,\eqno(3.33)$$
which implies $1\in I$ by the arguments in the above paragraph.
\psp

{\it Case 2}. $\hat{k}>0$ and there exists $p\in\{1,2\}$ such that $\vf_p(\al)=0$ and $j_p>0$ for some $a_{\al,\vec{j}}\neq 0$.
\psp

By (1.17), there exists $\tau\in \kn_{\vf_3}\bigcap \kn_{\vf_4}$ such that $\vf_p(\tau)=\vf_p(\al_0)$. 
We choose $\tau'\in\bigcap_{p'\neq q\in\ol{1,4}}\kn_{\vf_q}$ such that $\vf_{p'}(\tau'+\tau)\neq\vf_{p'}(\al_0),-\vf_{p'}(\al)$. Note that 
\begin{eqnarray*}& &({\cal A}_{4,\hat{k}-1}\bigcap I)\ni[x^{\tau'+\tau-\al_0,\vec{0}},u]_2\equiv\es_{p'}\vf_{p'}(\tau'+\tau-\al_0)\sum_{\vec{l}\in\vec{\cal J},|\vec{l}|=\hat{k}}l_pa_{\al,\vec{l}}x^{\al+\tau'+\tau,\vec{l}-1_{[p]}}\\&&+\sum_{\al+\tau'+\tau\neq \be\in \G,\vec{l}\in\vec{\cal J},|\vec{l}|=\hat{k}-1}c_{\be,\vec{l}}x^{\be,\vec{l}-1_{[p]}}\;\;(\mbox{mod}\:{\cal A}_{4,\hat{k}-2})\hspace{5.7cm}(3.34)\end{eqnarray*}
by (3.1) with $c_{\be,\vec{l}}\in \Bbb{F}$, which contradicts (3.5) because $\al+\tau'+\tau\neq \sgm$ due to $\vf_{p'}(\al+\tau'+\tau)\neq 0$.
\psp

{\it Case 3}. $\hat{k}>0$ and there exists $p\in\{3,4\}$ such that $\vf_p(\al+\al_0)=0$ and $j_p>0$ for some $a_{\al,\vec{j}}\neq 0$.
\psp

By (1.17), there exists $\tau\in \kn_{\vf_1}\bigcap \kn_{\vf_2}$ such that $\vf_p(\tau)=-\vf_p(\al_0)$. 
We choose $\tau'\in\bigcap_{p'\neq q\in\ol{1,4}}\kn_{\vf_q}\setminus\kn_{\vf_{p'}}$ such that $\vf_{p'}(\tau'+\tau+\al_0)\neq 0,-\vf_{p'}(\al)$. Note that 
\begin{eqnarray*}& &({\cal A}_{4,\hat{k}-1}\bigcap I)\ni[x^{\tau'+\tau,\vec{0}},u]_2\equiv\es_{p'}\vf_{p'}(\tau'+\tau+\al_0)\sum_{\vec{l}\in\vec{\cal J},|\vec{l}|=\hat{k}}l_pa_{\al,\vec{l}}x^{\al+\tau'+\tau,\vec{l}-1_{[p]}}\\&&+\sum_{\al+\tau'+\tau\neq \be\in \G,\vec{l}\in\vec{\cal J},|\vec{l}|=\hat{k}-1}c_{\be,\vec{l}}x^{\be,\vec{l}-1_{[p]}}\;\;(\mbox{mod}\:{\cal A}_{4,\hat{k}-2})\hspace{5.8cm}(3.35)\end{eqnarray*}
by (3.1) with $c_{\be,\vec{l}}\in \Bbb{F}$, which contradicts (3.5) because $\al+\tau'+\tau\neq \sgm$ due to $\vf_{p'}(\al+\tau'+\tau+\al_0)\neq 0$.

\psp

{\it Case 4}. $\hat{k}>0$. There  exist $p\in\{3,4\}$ and $a_{\al,\vec{j}}\neq 0$ such that $\vf_p(\al+\al_0)\neq 0$ and $j_p>0$.
\psp

We choose $\tau'\in\bigcap_{p'\neq q\in\ol{1,4}}\kn_{\vf_q}$ such that $\vf_{p'}(\tau'+\al_0)\neq 0$ and 
$$\vf_{p'}(\tau')\vf_p(\al+\al_0)+2\vf_{p'}(\al_0)\vf_p(\al)-2\vf_{p'}(\al)\vf_p(\al_0)\neq 0.\eqno(3.36)$$
Then we go back to Case 3 if $u$ is replaced by $[x^{\tau'-\al,\vec{0}},u]_2$, where
\begin{eqnarray*}& &[x^{\tau'-\al,\vec{0}},u]_2\equiv \es_{p'}[\vf_{p'}(\tau')\vf_p(\al+\al_0)+2\vf_{p'}(\al_0)\vf_p(\al)-2\vf_{p'}(\al)\vf_p(\al_0)]
\\& &\sum_{\vec{l}\in\vec{\cal J},|\vec{l}|=\hat{k}}a_{\al,\vec{l}}x^{\tau',\vec{l}}\;\;(\mbox{mod}\:{\cal A}_{\hat{k}-1})\hspace{9.2cm}(3.37)\end{eqnarray*}
by (3.1).
\psp

{\it Case 5}. $\hat{k}>0,\;\vf_3(\al_0)\vf_4(\al)-\vf_3(\al)\vf_4(\al_0)=0$ and there exist $p\in\{1,2\}$ and $a_{\al,\vec{j}}\neq 0$ such that $\vf_p(\al)\neq 0$ and $j_p>0$.
\psp

We choose $\tau'\in(\bigcap_{p'\neq q\in\ol{1,4}}\kn_{\vf_q})\setminus \kn_{\vf_{p'}}$ such that 
$$\vf_{p'}(\tau')\vf_p(\al)-\vf_{p'}(\al_0)\vf_p(\al)+\vf_{p'}(\al)\vf_p(\al_0)\neq 0.\eqno(3.38)$$
Then by  we go back to Case 2 if $u$ is replaced by $[x^{\tau'-\al-\al_0,\vec{0}},u]_2$, where
\begin{eqnarray*}& &[x^{\tau'-\al-\al_0,\vec{0}},u]_2\equiv \es_{p'}[\vf_{p'}(\tau')\vf_p(\al)-\vf_{p'}(\al_0)\vf_p(\al)+\vf_{p'}(\al)\vf_p(\al_0)]
\\& &\sum_{\vec{l}\in\vec{\cal J},|\vec{l}|=\hat{k}}a_{\al,\vec{l}}x^{\tau',\vec{l}}\;\;(\mbox{mod}\:{\cal A}_{\hat{k}-1})\hspace{9.2cm}(3.39)\end{eqnarray*}
by (3.1).

\psp

{\it Case 6}. $\hat{k}>0,\;\vf_3(\al_0)\vf_4(\al)-\vf_3(\al)\vf_4(\al_0)\neq 0$ and there exist $p\in\{1,2\}$ and $a_{\al,\vec{j}}\neq 0$ such that  $j_p>0$.
\psp

We go back to Case 2 if $u$ is replaced by $[x^{-\al,\vec{0}},u]_2$, where
$$[x^{-\al,\vec{0}},u]_2\equiv 2[\vf_3(\al_0)\vf_4(\al)-\vf_3(\al)\vf_4(\al_0)]
\sum_{\vec{l}\in\vec{\cal J},|\vec{l}|=\hat{k}}a_{0,\vec{l}}x^{0,\vec{l}}\;\;(\mbox{mod}\:{\cal A}_{\hat{k}-1})\eqno(3.40)$$
by (3.1).

\psp

This completes the proof of the conclusion in Step 1. 

\psp

{\it Step 2}. The conclusion of Step 1 implies $I={\cal A}$.
\psp

By (1.16), there exists $q\in\{3,4\}$ such that $\vf_q(\al_0)\neq 0$. Moreover, there exists $\tau_q\in \kn_{\vf_1}\bigcap\kn_{\vf_2}$ such that $\vf_q(\tau_q)=-\vf_q(\al_0)$ by (1.17). By (1.15), there exists $\tau'\in \bigcap_{q'\neq r\in\ol{1,4}}\kn_{\vf_r}\setminus\kn_{\vf_{q'}}$ such that $\vf_{q'}(\tau'+\tau_q+\al_0)\neq 0$, we have:
$$[x^{\tau'+\tau_q,\vec{0}},1]_2=\es_{q'}\vf_{q'}(\tau'+\tau_q+\al_0)\vf_q(\al_0)x^{\tau'+\tau_q,\vec{0}}\in I\eqno(3.41)$$
by (3.1). So $x^{\tau'+\tau_q,\vec{0}}\in I$. Moreover, for any $(\al,\vec{i})\in \G\times\vec{\cal J}$, we have
$$[x^{\tau'+\tau_q,\vec{0}},x^{\al,\vec{i}}]_2=\es_{q'}\vf_{q'}(\tau'+\tau_q+\al_0)(\vf_q(\al+\al_0)x^{\al+\tau_q+\tau',\vec{i}}+i_qx^{\al+\tau_q+\tau',\vec{i}-1_{[q]}}).\eqno(3.42)$$
By (3.42) and induction on $i_q$, we can prove $I={\cal A}_4$ if ${\cal J}_q=\Bbb{N}$. Assume
that ${\cal J}_q=\{0\}$. Then (3.42) shows
$$x^{\be,\vec{j}}\in I\qquad\mbox{for any}\;\;(\be,\vec{j})\in \G\times\vec{\cal J},\;\vf_q(\be)\neq -2\vf_q(\al_0)\eqno(3.43)$$
because $\vf_q(\al+\tau_q+\tau')=-2\vf_q(\al_0)$ if $\vf_q(\al+\al_0)=0$. 
Furthemore, there exists $\tau_{q'}\in \kn_{\vf_1}\bigcap\kn_{\vf_2}$ such that $\vf_{q'}(\tau_{q'})=-\vf_{q'}(\al_0)$ by (1.17). Thus $x^{\tau+\tau_{q'}}\in I$ for some $\tau\in \bigcap_{q\neq r\in\ol{1,4}}\kn_{\vf_r}\setminus\kn_{\vf_p}$ such that $\vf_q(\tau+\tau_{q'})\neq -\vf_q(\al_0),-2\vf_q(\al_0)$ by (1.15). Exchanging positions of $q$ and $q'$ in (3.42), we can prove $I={\cal A}_4$ if ${\cal J}_{q'}=\Bbb{N}$ and (3.43) with $q$ replaced by $q'$ if ${\cal J}_{q'}=\{0\}$.
Thus we always have $I={\cal A}_4$ if ${\cal J}_3=\Bbb{N}$ or ${\cal J}_4=\Bbb{N}$.

Assume that ${\cal J}_3={\cal J}_4=\{0\}$. The above paragraph has proved
$$x^{\be,\vec{j}}\in I\qquad\mbox{for any}\;\;(\be,\vec{j})\in \G\times\vec{\cal J},\;\vf_3(\be)\neq -2\vf_3(\al_0)\;\mbox{or}\;\vf_4(\be)\neq -2\vf_4(\al_0).\eqno(3.44)$$
Let $p\in\{1,2\}$. By (1.17), there exists $\tau_p\in \kn_{\vf_3}\bigcap\kn_{\vf_4}$ such that $\vf_p(\tau_p)=\vf_p(\al_0)$. Choose $\tau'\in \bigcap_{p'\neq r\in\ol{1,4}}\kn_{\vf_r}\setminus \kn_{\vf_{p'}}$ such that $\vf_{p'}(\tau'+\tau_p-\al_0)\neq 0$ by (1.15), we have $x^{\tau'+\tau_p-\al_0,\vec{0}}\in I$ by (3.44) and
$$[x^{\tau'+\tau_p-\al_0,\vec{0}},x^{\al,\vec{i}}]=\es_{p'}\vf_{p'}(\tau'+\tau_p-\al_0)(\vf_p(\al)x^{\al+\tau_p+\tau',\vec{i}}+i_px^{\al+\tau_p+\tau',\vec{i}-1_{[q]}})\eqno(3.45)$$
for any $(\al,\vec{i})\in \G\times\vec{\cal J}$.  By (3.45) and induction on $i_p$, we can prove $I={\cal A}_4$ if ${\cal J}_p=\Bbb{N}$. Assume that ${\cal J}_p=\{0\}$. Then (3.45) shows 
$$x^{\be,\vec{j}}\in I\qquad\mbox{for any}\;\;(\be,\vec{j})\in \G\times\vec{\cal J},\;\vf_p(\be)\neq \vf_p(\al_0)\eqno(3.46)$$
because $\vf_p(\al+\tau_p+\tau')=\vf_p(\tau_p)=\vf_p(\al_0)$ if $\vf_p(\al)=0$. Since (1.20) is equivalent to
$$\vf_1(\rho)=\vf_1(\al_0),\;\vf_2(\rho)=\vf_2(\al_0),\;\vf_3(\rho)=-2\vf_3(\al_0),\;\vf_4(\rho)=-2\vf_4(\al_0),\eqno(3.47)$$
We obtain the first Statement in Theorem 2.
\psp

{\it Proof of the Second Statement in Theorem}
\psp

Now $\vec{\cal J}=\{\vec{0}\}$ and there exists $\rho\in \G$ such that (1.20) holds. 
Set
$$\hat{\cal B}_4=\sum_{\rho\neq \al\in\G}\Bbb{F}x^{\al,\vec{0}}.\eqno(3.48)$$
For any $\al,\be\in\G$, we have:
\begin{eqnarray*}[x^{\al,\vec{0}},x^{\be,\vec{0}}]_2&=&(\vf_1(\al)\vf_2(\be)-\vf_1(\be)\vf_2(\al))x^{\al+\be+\al_0,\vec{0}}+(\vf_3(\al+\al_0)\vf_4(\be+\al_0)\\& &-\vf_3(\be+\al_0)\vf_4(\al+\al_0))x^{\al+\be,\vec{0}}\hspace{6.4cm}(3.49)\end{eqnarray*}
by (3.1). If $\vf_1(\al+\be+\al_0)=\vf_1(\al_0)$ and $\vf_2(\al+\be+\al_0)=\vf_2(\al_0)$, then $\vf_1(\al)=-\vf_1(\be)$ and $\vf_2(\al)=-\vf_2(\be)$. Thus $\vf_1(\al)\vf_2(\be)-\vf_1(\be)\vf_2(\al)=0$. If $\vf_3(\al+\be)=-2\vf_3(\al_0)$ and $\vf_4(\al+\be)=-2\vf_4(\al_0)$, then $\vf_3(\al)=-\vf_3(\be)-2\vf_3(\al_0)$ and $\vf_4(\al)=-\vf_4(\be)-2\vf_4(\al_0)$. So
\begin{eqnarray*}& &\vf_3(\al+\al_0)\vf_4(\be+\al_0)-\vf_3(\be+\al_0)\vf_4(\al+\al_0)\\&=&(\vf_3(\al)+\vf_3(\al_0))\vf_4(\be+\al_0)-\vf_3(\be+\al_0)(\vf_4(\al)+\vf_4(\al_0))\\&=&(-\vf_3(\be)-\vf_3(\al_0))(\vf_4(\be)+\vf_4(\al_0))-(\vf_3(\be)+\vf_3(\al_0))(-\vf_4(\be)-\vf_4(\al_0))\\&=&0.\hspace{13.6cm}(3.50)\end{eqnarray*}
Thus $[{\cal A}_4,{\cal A}_4]_2\subset \hat{\cal B}_4$. Hence $[{\cal A}_4,{\cal A}_4]_2=\hat{\cal B}_4$ by (3.42) and (3.45) with $\vec{i}=\vec{0}$. Replacing ${\cal A}_4$ by $\hat{\cal B}_4$ and $\G$ by $\G\setminus\{\rho\}$ in the proof of the first statement in Thorem 2, we obtain (3.41) and (3.43), which implies $I=\hat{\cal B}_4$. Therefore, 
$${\cal B}_4^{(1)}=[{\cal B}_4,{\cal B}_4]_2=\hat{\cal B}_4/\Bbb{F}x^{\sgm,\vec{0}}\eqno(3.51)$$
is simple. 

This completes the proof of Theorem 2.

\subsection{Related Examples of Simple Lie Algebras}

Let $\Bbb{R}$ be the field of real numbers. By Theorem 2, we obtain the following examples of simple Lie algebras $({\cal B},[\cdot,\cdot])$, where $f_t$ denotes the partial derivative of a polynomial $f$ with respect to the variable $t$. Let $m$ be a positive integer and let $n$ be a nonzero integer.
\psp

{\bf Example 3.1}. Let ${\cal A}$ be the subspace
$$\sum_{j=1}^{m-1}(t_1t_2t_3t_4)^{j/m}\Bbb{R}[t_i,t_i^{-1}\mid i\in\ol{1,4}]+\sum_{(n,n,-2n,-2n)\neq (l_1,l_2,l_3,l_4)\in\Bbb{Z}^4}\Bbb{R}t_1^{l_1}t_2^{l_2}t_3^{l_3}t_4^{l_4}\eqno(3.52)$$
of the polynomial algebra $\Bbb{R}[t_i^{1/m},t_i^{-1/m}\mid i\in\ol{1,4}]$. The space ${\cal B}={\cal A}/\Bbb{R}(t_3t_4)^{-n}$ and its Lie bracket is induced by
\begin{eqnarray*}[f,g]&=&t_1t_2(t_1t_2t_3t_4)^n(f_{t_1}g_{t_2}-f_{t_2}g_{t_1})+(t_3f_{t_3}+nf)(t_4g_{t_4}+ng)\\& &-(t_4f_{t_4}+nf)(t_3g_{t_3}+ng)\hspace{8.5cm}(3.53)\end{eqnarray*}
for $f,g\in{\cal A}$.
\pse

{\bf Example 3.2}. Let ${\cal A}$ be the subalgebra of 
$$\Bbb{R}[t_i^{1/m},t_i^{-1/m},t_5\mid i=\ol{1,4}]\eqno(3.54)$$
 generated by 
$$\{t_i,t_i^{-1},t_5,(t_1t_2t_3t_4)^{1/m}\mid i=\ol{1,4}\}.\eqno(3.55)$$
When $m=1$, 
$${\cal A}=\Bbb{R}[t_i,t_i^{-1},t_5\mid i=\ol{1,4}].\eqno(3.56)$$
The space ${\cal B}={\cal A}/\Bbb{R}(t_3t_4)^{-n}$ and its Lie bracket is induced by
\begin{eqnarray*}[f,g]&=&t_1t_2(t_1t_2t_3t_4)^n(f_{t_1}g_{t_2}-f_{t_2}g_{t_1})+(t_3f_{t_3}+f_{t_5}+nf)(t_4g_{t_4}+ng)\\& &-(t_4f_{t_4}+nf)(t_3g_{t_3}+g_{t_5}+ng)\hspace{7.5cm}(3.57)\end{eqnarray*}
for $f,g\in{\cal A}$.
\pse

{\bf Example 3.3}. Let ${\cal A}$ be the subalgebra of 
$$\Bbb{R}[t_i^{1/m},t_i^{-1/m},t_5,t_6\mid i=\ol{1,4}]\eqno(3.58)$$
 generated by 
$$\{t_i,t_i^{-1},t_5,t_6,(t_1t_2t_3t_4)^{1/m}\mid i=\ol{1,4}\}.\eqno(3.59)$$
When $m=1$,
$${\cal A}=\Bbb{R}[t_i,t_i^{-1},t_5,t_6\mid i=\ol{1,4}].\eqno(3.60)$$
The space ${\cal B}={\cal A}/\Bbb{R}(t_3t_4)^{-n}$ and its Lie bracket is induced by
\begin{eqnarray*}[f,g]&=&t_1t_2(t_1t_2t_3t_4)^n(f_{t_1}g_{t_2}-f_{t_2}g_{t_1})+(t_3f_{t_3}+f_{t_5}+nf)(t_4g_{t_4}+g_{t_6}+ng)\\& &-(t_4f_{t_4}+f_{t_6}+nf)(t_3g_{t_3}+g_{t_5}+ng)\hspace{6.6cm}(3.61)\end{eqnarray*}
or
\begin{eqnarray*}[f,g]&=&t_2(t_1t_2t_3t_4)^n[(t_1f_{t_1}+f_{t_5})g_{t_2}-f_{t_2}(t_1g_{t_1}+g_{t_5})]\\& &+(t_3f_{t_3}+f_{t_6}+nf)(t_4g_{t_4}+ng)-(t_4f_{t_4}+nf)(t_3g_{t_3}+g_{t_6}+ng)\hspace{1.9cm}(3.62)\end{eqnarray*}
for $f,g\in{\cal A}$.
\pse

{\bf Example 3.4}. Let ${\cal A}$ be the subalgebra of 
$$\Bbb{R}[t_i^{1/m},t_i^{-1/m},t_5,t_6,t_7\mid i=\ol{1,4}]\eqno(3.63)$$
 generated by 
$$\{t_i,t_i^{-1},t_5,t_6,t_7,(t_1t_2t_3t_4)^{1/m}\mid i=\ol{1,4}\}.\eqno(3.64)$$
When $m=1$,
$${\cal A}=\Bbb{R}[t_i,t_i^{-1},t_5,t_6,t_7\mid i=\ol{1,4}].\eqno(3.65)$$
The space ${\cal B}={\cal A}/\Bbb{R}(t_3t_4)^{-n}$ and its Lie bracket is induced by
\begin{eqnarray*}[f,g]&=&t_2(t_1t_2t_3t_4)^n[(t_1f_{t_1}+f_{t_5})g_{t_2}-f_{t_2}(t_1g_{t_1}+g_{t_5})]+(t_3f_{t_3}+f_{t_6}+nf)\\&&\times(t_4g_{t_4}+g_{t_7}+ng)-(t_4f_{t_4}+f_{t_7}+nf)(t_3g_{t_3}+g_{t_6}+ng)\hspace{3cm}(3.66)\end{eqnarray*}
for $f,g\in{\cal A}$.
\pse

{\bf Example 3.5}. Let ${\cal A}$ be the subalgebra of 
$$\Bbb{R}[t_i^{1/m},t_i^{-1/m},t_{4+i}\mid i=\ol{1,4}]\eqno(3.67)$$
 generated by 
$$\{t_i,t_i^{-1},t_{4+i},(t_1t_2t_3t_4)^{1/m}\mid i=\ol{1,4}\}.\eqno(3.68)$$
When $m=1$,
$${\cal A}=\Bbb{R}[t_i,t_i^{-1},t_5,t_6,t_7\mid i=\ol{1,4}].\eqno(3.69)$$
The space ${\cal B}={\cal A}/\Bbb{R}(t_3t_4)^{-n}$ and its Lie bracket is induced by
\begin{eqnarray*}& &[f,g]=\\& &(t_1t_2t_3t_4)^n[(t_1f_{t_1}+f_{t_5})(t_2g_{t_2}+g_{t_6})-(t_2f_{t_2}+f_{t_6})(t_1g_{t_1}+g_{t_5})]+(t_3f_{t_3}+f_{t_7}+nf)\\& &\times(t_4g_{t_4}+g_{t_8}+ng)-(t_4f_{t_4}+f_{t_8}+nf)(t_3g_{t_3}+g_{t_7}+ng)\hspace{4.2cm}(3.70)\end{eqnarray*}
for $f,g\in{\cal A}$.

\section{Proof of Theorem 3 and Examples}

In this section, we shall first give the proof of Theorem 3 and then present some related examples of simple Lie superalgebras.

\subsection{Proof of Theorem 3}

We assume that the conditions  in Theorem 3 hold. First we have
\begin{eqnarray*}& &[(x^{\al,\vec{i}})_{[0]},(x^{\be,\vec{j}})_{[1]}]\\&=&(\vf_1(\al)\vf_2(\be)-\vf_1(\be)\vf_2(\al)+\ves\vf_1(\be)+\ves \vf_1(\al_0-\al)/2+(\vf_2(\al_0)\vf_1(\al)\\& &-\vf_1(\al_0)\vf_2(\al))/2)
(x^{\al+\be,\vec{i}+\vec{j}})_{[1]}+(i_1\vf_2(\be)-j_1\vf_2(\al)+\ves(j_1-i_1/2)\\& &+i_1\vf_2(\al_0)/2)(x^{\al+\be,\vec{i}+\vec{j}-1_{[1]}})_{[1]}+(i_1j_2-i_2j_1)(x^{\al+\be,\vec{i}+\vec{j}-1_{[1]}-1_{[2]}})_{[1]}\\& &+(j_2\vf_1(\al)-i_2\vf_1(\be)-i_2\vf_1(\al_0)/2)(x^{\al+\be,\vec{i}+\vec{j}-1_{[2]}})_{[1]}\hspace{5.1cm}(4.1)\end{eqnarray*}
for $(\al,\vec{i}),(\be,\vec{j})\in \G\times\vec{\cal J}$ by (1.26).
 In particular,
\begin{eqnarray*}[(x^{\al,\vec{0}})_{[0]},(x^{\be,\vec{j}})_{[1]}]&=&(\vf_1(\al)\vf_2(\be)-\vf_1(\be)\vf_2(\al)+\ves\vf_1(\be)+\ves \vf_1(\al_0-\al)/2\\& &+(\vf_2(\al_0)\vf_1(\al)-\vf_1(\al_0)\vf_2(\al))/2)(x^{\al+\be,\vec{j}})_{[1]}\\& &+j_1(\ves-\vf_2(\al))(x^{\al+\be,\vec{j}-1_{[1]}})_{[1]}+j_2\vf_1(\al)(x^{\al+\be,\vec{j}-1_{[2]}})_{[1]}\hspace{1.6cm}(4.2)\end{eqnarray*}
for $\al\in\G$ and $(\be,\vec{j})\in \G\times\vec{\cal J}$. When $\al=0$, we get
$$[1_{[0]},(x^{\be,\vec{j}})_{[1]}]_2=\ves(\vf_1(\be)+\vf_1(\al_0)/2)(x^{\be,\vec{j}})_{[1]}+\ves j_1(x^{\be,\vec{j}-1_{[1]}})_{[1]}.\eqno(4.3)$$

If $\ves=1$ and ${\cal J}_1=\Bbb{N}$, we get
$$\td{\cal B}_1=[\td{\cal A}_0,\td{\cal A}_1]=\td{\cal A}_1\eqno(4.4)$$
by (4.3) and induction on $j_1$. 
 Taking 
$\al\in\kn_{\vf_2}\setminus \kn_{\vf_1}$ in (4.2) when $\ves =0$ and $\vf_1\not\equiv 0$ by (1.29), we get
\begin{eqnarray*}[(x^{\al,\vec{0}})_{[0]},(x^{\be,\vec{j}})_{[1]}]&=&\vf_1(\al)(\vf_2(\be)+\vf_2(\al_0)/2)
(x^{\al+\be,\vec{j}})_{[1]}\\& &+j_2\vf_1(\al)(x^{\al+\be,\vec{j}-1_{[2]}})_{[1]}\hspace{7cm}(4.5)\end{eqnarray*}
for $(\be,\vec{j})\in \G\times\vec{\cal J}$. If ${\cal J}_2=\Bbb{N}$, then we obtain (4.4) by (4.2) when ${\cal J}_1=\{0\}$ and by (4.5) when $\ves=0$ and $\vf_1\not\equiv 0$, and by induction on $j_2$. Consider the case that $\ves=0,\;\vf_2\not\equiv 0$ and ${\cal J}_1=\Bbb{N}$. Letting $\al\in\kn_{\vf_1}\setminus \kn_{\vf_2}$ in (4.2) by (1.29), we have
$$[(x^{\al,\vec{0}})_{[0]},(x^{\be,\vec{j}})_{[1]}]=-\vf_2(\al)[(\vf_1(\be)+\vf_1(\al_0)/2)(x^{\al+\be,\vec{j}})_{[1]}+j_1(x^{\al+\be,\vec{j}-1_{[1]}})_{[1]}]\eqno(4.6)$$
for $(\be,\vec{j})\in \G\times\vec{\cal J}$. Again we have (4.4) by (4.6) and induction on $j_1$. Suppose that $\ves=0,\;\vf_1\equiv\vf_2\equiv 0$ and ${\cal J}_1={\cal J}_2=\Bbb{N}$, we have
$$[(x^{0,1_{[1]}})_{[0]},(x^{0,\vec{j}})_{[1]}]=j_2(x^{0,\vec{j}-1_{[2]}})_{[1]}\eqno(4.7)$$
for $\vec{j}\in\vec{J}$ by (4.1), which implies (4.4). Thus we have proved (4.4) holds when $\vec{\cal J}\neq\{\vec{0}\}$.

Assume that $\vec{\cal J}=\{\vec{0}\}$. By (4.2),
\begin{eqnarray*}& &[(x^{\al,\vec{0}})_{[0]},(x^{\be,\vec{0}})_{[1]]}]=[\vf_1(\al)\vf_2(\be)-\vf_1(\be)\vf_2(\al)+\ves\vf_1(\be)+\ves \vf_1(\al_0-\al)/2\\& &+(\vf_2(\al_0)\vf_1(\al)-\vf_1(\al_0)\vf_2(\al))/2](x^{\al+\be,\vec{0}})_{[1]}\hspace{6.5cm}(4.8)\end{eqnarray*}
for $\al,\be\in\G$. In particular,
$$[1_{[0]},(x^{\be,\vec{0}})_{[1]}]=(\ves\vf_1(\be)+\ves \vf_1(\al_0)/2)(x^{\be,\vec{0}})_{[1]},\eqno(4.9)$$
$$[(x^{\al,\vec{0}})_{[0]},1_{[1]}]={1\over 2}(\ves \vf_1(\al_0-\al)+(\vf_2(\al_0)\vf_1(\al)-\vf_1(\al_0)\vf_2(\al)))(x^{\al,\vec{0}})_{[1]}.\eqno(4.10)$$
Suppose that $\ves=0,\;\al_0\neq 0$ and (1.29) holds. By (4.10),
$$(x^{\al,\vec{0}})_{[1]}\in \td{\cal B}_1\qquad\mbox{for}\;\;\al\in\G,\;\vf_2(\al_0)\vf_1(\al)-\vf_1(\al_0)\vf_2(\al)\neq 0.\eqno(4.11)$$
Note that $\vf_1\not\equiv 0$ and $\vf_2\not\equiv 0$ by (1.6). 
For any $\be,\tau\in\G$ such that $\vf_2(\al_0)\vf_1(\be)-\vf_1(\al_0)\vf_2(\be)=0$, we have
$$[(x^{-\tau+\be,\vec{0}}_{[0]},(x^{\tau,\vec{0}})_{[1]}]=[(\vf_1(\be)+\vf_1(\al_0)/2)\vf_2(\tau)-(\vf_2(\be)+\vf_2(\al_0)/2)\vf_1(\tau)](x^{\be,\vec{0}})_{[1]}.\eqno(4.12)$$
If 
$$\vf_1(\be)+\vf_1(\al_0)/2\neq 0\;\;\mbox{or}\;\;\vf_2(\be)+\vf_2(\al_0)/2\neq 0,\eqno(4.13)$$
then we can choose $\tau\in \G$ such that
$$(\vf_1(\be)+\vf_1(\al_0)/2)\vf_2(\tau)-(\vf_2(\be)+\vf_2(\al_0)/2)\vf_1(\tau)\neq 0\eqno(4.14)$$
by (1.29) and have $(x^{\be,\vec{0}})_{[1]}\in \td{\cal B}_1$. Thus we have
$$\td{\cal B}_1=\{(x^{\al,\vec{0}})_{[1]}\mid \al\in\G,\;2\al\neq -\al_0\}.\eqno(4.15)$$

Next we assume that $\ves=1$ and $\vf_1(\al_0)=0$. Then by (4.9),
$$(x^{\al,\vec{0}})_{[1]}\in \td{\cal B}_1\qquad\mbox{for}\;\;\al\in\G,\;\vf_1(\al)\neq 0.\eqno(4.16)$$
For any  $\be\in\kn_{\vf_1}$ and $\tau\in\G$, we have
$$[(x^{-\tau+\be,\vec{0}})_{[0]},(x^{\tau,\vec{0}})_{[1]}]={1\over 2}\vf_1(\tau)(3-\vf_2(\al_0)-2\vf_2(\be))(x^{\be,\vec{0}})_{[1]}\eqno(4.17)$$
by (4.2). If $3-\vf_2(\al_0)-2\vf_2(\be)\neq 0$, then we choose $\tau\in\kn_{\vf_2}\setminus \kn_{\vf_1}$ in (4.17) by (1.29) and have $(x^{\be,\vec{0}})_{[1]}\in \td{\cal B}_1$. Thus we have
$$\td{\cal B}_1=\{(x^{\al,\vec{0}})_{[1]}\mid \al\in\G,\;\vf_1(\al)\neq 0\;\mbox{or}\;2\vf_2(\al)\neq 3-\vf_2(\al_0)\}.\eqno(4.18)$$

Now we assume that $\ves=1$ and $\vf_1(\al_0)\neq 0$. By (4.9) and (4.10),
$$(x^{\al,\vec{0}})_1\in \td{\cal B}_1\;\;\mbox{for}\;\;\al\in\G,\;\vf_1(\al)\neq -{\vf_1(\al_0)\over 2}\;\;\mbox{or}\;\;\vf_2(\al)\neq {3-\vf_2(\al_0)\over 2}.\eqno(4.19)$$
By (1.7), there exists at most one $\kappa\in\G$ such that
$$\vf_1(\kappa)= -{\vf_1(\al_0)\over 2},\qquad \vf_2(\kappa) = {3-\vf_2(\al_0)\over 2}.\eqno(4.20)$$
Assume that such $\kappa$ exists. By (4.8),
\begin{eqnarray*}& &[(x^{\al,\vec{0}})_{[0]},(x^{\kappa-\al,\vec{0}})_{[1]]}]\\&=&[\vf_1(\al)\vf_2(\kappa)-\vf_1(\kappa)\vf_2(\al)+\vf_1(\kappa-\al)+\vf_1(\al_0-\al)/2\\& &+(\vf_2(\al_0)\vf_1(\al)-\vf_1(\al_0)\vf_2(\al))/2](x^{\kappa,\vec{0}})_{[1]}\\&=& (\vf_1(\al)(3-\vf_2(\al_0))/2+\vf_1(\al_0)\vf_2(\al)/2-\vf_1(\al_0)/2-\vf_1(\al)+\vf_1(\al_0-\al)/2\\& &+(\vf_2(\al_0)\vf_1(\al)-\vf_1(\al_0)\vf_2(\al))/2)(x^{\kappa,\vec{0}})_{[1]}\\&=&0\hspace{13.7cm}(4.21)\end{eqnarray*}
for any $\al\in\G$. Thus 
$$\td{\cal B}_1=\{(x^{\al,\vec{0}})_{[1]}\mid \al\in\G,\;\al\neq \kappa\}.\eqno(4.22)$$

Therefore, the codimension of $\td{\cal B}$ in $\td{\cal A}$ is at most one. 

By (1.26), $1_{[0]}$ is a central element of $\td{\cal A}$ when $\ves=0$. Recall that $\sgm_1$ is an element of $\kn_{\vf_1}$ satisfying $\vf_2(\sgm_1)=1$. When $\ves=1$, we have
$$[(x^{\sgm_1,\vec{0}})_{[0]},v_{[1]}]=(1-\vf_2(\sgm_1))(x^{\sgm_1,\vec{0}}\ptl_1(v))_{[1]}+{\vf_1(\al_0)(1-\vf_2(\sgm_1))\over 2}(x^{\sgm_1,\vec{0}}v)_{[1]}=0\eqno(4.23)$$
by (1.12) and (1.26). Expression (4.23) shows that $(x^{\sgm_1,\vec{0}})_{[0]}$ is a central
element of $\td{\cal A}$ when $\ves=1$. Thus (1.28) defines a quotient algebra $\td{\cal C}$. 

Let $I$ be any ideal of $\td{\cal B}$ such that $I$ strictly contains $\Bbb{F}1_{[0]}$ when $\ves=0$ and $\Bbb{F}(x^{\sgm_1,\vec{0}})_{[0]}$ when $\ves=1$. To prove Theorem 3 is equivalent to proving $I=\td{\cal B}$.

First we assume that $I\bigcap \td{\cal B}_1\neq\{0\}$. Note that
$$I\supset [I\bigcap \td{\cal B}_1,\td{\cal B}_1]\not\subset \Bbb{F}+\Bbb{F}x^{\sgm_1,\vec{0}}+\Bbb{F}x^{\sgm_2,\vec{0}}\eqno(4.24)$$
by (1.25), (4.4), (4.15), (4.18) and (4.22). Moreover, it is known that
$$\td{A}_0\;\;\mbox{is a simple Lie algebra when}\;\;\vf_2\equiv 0,\;{\cal J}_2=\{0\}\;\mbox{and}\;\ves=1\eqno(4.25)$$
(e.g., cf. Theorem 2.2 in [X4]) and 
$$\td{A}_0/\Bbb{F}1_{[0]}\;\;\mbox{forms a simple Lie algebra}\eqno(4.26)$$
when $\ves=0$ by theorem 4.1 in [X4]. Thus
$$\td{\cal A}_0\subset I \qquad \mbox{if}\;\;\ves=0\eqno(4.27)$$
by (4.24) and (4.26). Moreover, by Section 2, (2.21) and (4.25), we have
$$(x^{\al,\vec{i}})_{[0]}\in I\;\;\mbox{for}\;\;(\sgm_2,\vec{0})\neq (\al,\vec{i})\in \G\times\vec{\cal J}\eqno(4.28)$$
if $\ves=1$, where $\sgm_2\in\kn_{\vf_1}$ and $\vf_2(\sgm_2)=2$.  Furthermore, by the arguments in (4.1)-(4.22), we have
$$\td{\cal B}_1\subset I.\eqno(4.26)$$
Note that there exists $l\in\Bbb{N}$ such that $(x^{(l+1)\sgm_2,\vec{0}})_{[1]}\in\td{\cal B}_1$ by (4.4), (4.15), (4.18) and (4.22). Hence, we have
$$[(x^{-l\sgm_2-\al_0,\vec{0}})_{[1]},(x^{(l+1)\sgm_2,\vec{0}})_{[1]}]=x^{\sgm_2,\vec{0}}\in I\eqno(4.30)$$
by (1.25). Thus $I=\td{\cal B}$.

Next we assume
$$I\bigcap \td{\cal A}_0\not\in \Bbb{F}1_{[0]}\;\;\mbox{if}\;\;\ves=0\eqno(4.31)$$
and 
$$I\bigcap \td{\cal A}_0\not\in \Bbb{F}(x^{\sgm_1,\vec{0}})_{[0]}\;\;\mbox{if}\;\;\ves=1.\eqno(4.32)$$
Then we have (4.27) and (4.28) by Section 2, (4.25) and (2.26), which implies $I=\td{\cal B}$.

Now we assume that $I\bigcap \td{\cal B}_1=\{0\}$ and $I\not\subset \td{\cal A}_0$. Let $w=u_{[0]}+v_{[1]}\in I\setminus \td{\cal A}_0$ with $u,v\in {\cal A}_2$. Then $v\neq 0$ and $u\not\in\Bbb{F}1_{[0]}$ when $\ves=0$ and $u\not\in\Bbb{F}(x^{\sgm_1,\vec{0}})_{[0]}$ when $\ves=1$. 
Since
$$[u_{[0]}, w]=[u_{[0]},v_{[1]}]\in I\bigcap \td{\cal B}_1,\eqno(4.33)$$
we have $[u_{[0]},v_{[1]}]=0$. Thus
$$[v_{[1]}, w]=[v_{[1]},v_{[1]}]=(x^{\al_0,\vec{0}}v^2)_{[0]}\in I.\eqno(4.34)$$
If $x^{\al_0,\vec{0}}v^2\not\in\Bbb{F}$ when $\ves=0$ and $x^{\al_0,\vec{0}}v^2\not\in\Bbb{F}x^{\sgm_1,\vec{0}}$ when $\ves=1$ (which naturally holds if $\vf_2\equiv 0$), then we have (4.27) and (4.28) by Section 2, (4.25) and (4.26), which implies $I=\td{\cal B}$ that contradicts $I\bigcap \td{\cal B}_1=\{0\}$.

Consider the case that $\ves=0$ and $x^{\al_0,\vec{0}}v^2\in\Bbb{F}$. Up to a constant multiple, we can assume that
$$v=x^{-\al_0/2,\vec{0}}.\eqno(4.35)$$
Note that by (4.1),
$$[(x^{\al,\vec{i}})_{[0]},(x^{-\al_0/2,\vec{0}})_{[1]}]=0\qquad\mbox{for any}\;\;(\vec{\al},\vec{i})\in \G\times \vec{\cal J}.\eqno(4.36)$$
Thus
$$[\td{\cal A}_0,w]\in I\bigcap \td{\cal A}_0,\;\;[\td{\cal A}_0,w]\not\subset\Bbb{F}1_{[0]}\eqno(4.37)$$
 by (4.26), which leads to $I=\td{\cal B}$ that contradicts $I\bigcap \td{\cal B}_1=\{0\}$.

Consider the case that $\ves=1,\;\vf_2\not\equiv 0$ and $x^{\al_0,\vec{0}}v^2\in\Bbb{F}x^{\sgm_1,\vec{0}}$. Up to a constant multiple, we can assume that
$$v=x^{(\sgm_1-\al_0)/2,\vec{0}}.\eqno(4.38)$$
Note that
$$[(x^{\al,\vec{i}})_{[0]},(x^{(\sgm_1-\al_0)/2,\vec{0}})_{[1]}]=0\qquad\mbox{for any}\;\;(\vec{\al},\vec{i})\in \G\times \vec{\cal J}\eqno(4.39)$$
by (4.1). So we have 
$$[\td{\cal A}_0,w]\in I\bigcap \td{\cal A}_0,\;\;[\td{\cal A}_0,w]\not\subset\Bbb{F}(x^{\sgm_1,\vec{0}})_{[0]}\eqno(4.40)$$
by Section 2, which leads to $I=\td{\cal B}$ that contradicts $I\bigcap \td{\cal B}_1=\{0\}$. This completes the proof of Theorem 3.

\subsection{Related Examples of Simple Lie Superalgebras}

Let $\Bbb{R}$ be the field of real numbers. By Theorem 3, we obtain the following examples of simple Lie superalgebras $(\bar{\cal B},[\cdot,\cdot])$ (cf. (1.21) and (1.22)).
\psp

{\bf Example 4.1}. The space $\bar{\cal B}=\Bbb{R}[t]\times \Bbb{R}[t]$ with the Lie bracket:
$$[f_{[0]},g_{[0]}]=(fg'-f'g)_{[0]},\;\;[f_{[1]},g_{[1]}]=(fg)_{[0]},\;\;[f_{[0]},g_{[1]}]=(fg'-f'g/2)_{[1]}\eqno(4.41)$$
for $f,g\in \Bbb{R}[t]$ (cf. (1.23)), where $f'=df/dt$.
\pse

{\bf Example 4.2}. The space $\bar{\cal B}=\Bbb{R}[t,t^{-1}]\times \Bbb{R}[t,t^{-1}]$ with the Lie bracket:
$$[f_{[0]},g_{[0]}]=(tfg'-tf'g)_{[0]},\;[f_{[1]},g_{[1]}]=(tfg)_{[0]},\;[f_{[0]},g_{[1]}]=[tfg'+(f-tf')g/2]_{[1]}\eqno(4.42)$$
for $f,g\in \Bbb{R}[t,t^{-1}]$ (cf. (1.23)). In this case, the Lie superalgebra $(\bar{\cal B},[\cdot,\cdot])$ is the centerless super Virasoro algebra.

{\bf Example 4.3}. The space $\bar{\cal B}=\Bbb{R}[t_1,t_2]\times \Bbb{R}[t_1,t_2]$ with the Lie bracket:
$$[f_{[0]},g_{[0]}]=[f_{t_1}(g_{t_2}-g)+(f-f_{t_2})g_{t_1}]_{[0]},\;\;[f_{[1]},g_{[1]}]=(fg)_{[0]},\eqno(4.43)$$
$$[f_{[0]},g_{[1]}]=[f_{t_1}(g_{t_2}-g/2)+(f-f_{t_2})g_{t_1}]_{[1]}\eqno(4.44)$$
for $f,g\in \Bbb{R}[t_1,t_2,t_1^{-1}]$ (cf.  (1.23)). 
\pse

{\bf Example 4.4}. Let $n$ be an integer. The space $\bar{\cal B}=\Bbb{R}[t_1,t_2,t_1^{-1}]\times \Bbb{R}[t_1,t_2,t_1^{-1}]$ with the Lie bracket:
$$[f_{[0]},g_{[0]}]=[(t_1g_{t_1}+g_{t_2})f-(t_1f_{t_1}+f_{t_2})g]_{[0]},\;\;[f_{[1]},g_{[1]}]=(t_1^nfg)_{[0]},\eqno(4.45)$$
$$[f_{[0]},g_{[1]}]=[(t_1g_{t_1}+g_{t_2})f+(nf-t_1f_{t_1}-f_{t_2})g/2]_{[1]}\eqno(4.46)$$
or
$$[f_{[0]},g_{[0]}]=[t_1f_{t_1}(g_{t_2}-g)+t_1(f-f_{t_2})g_{t_1}]_{[0]},\;\;[f_{[1]},g_{[1]}]=(t_1^nfg)_{[0]},\eqno(4.47)$$
$$[f_{[0]},g_{[1]}]=[t_1f_{t_1}(g_{t_2}-g/2)+(f-f_{t_2})(t_1g_{t_1}+ng/2)]_{[1]}\eqno(4.48)$$
for $f,g\in \Bbb{R}[t_1,t_2,t_1^{-1}]$ (cf. (1.23)). 
The space 
$$\bar{\cal B}=(\Bbb{R}[t_1,t_2,t_1^{-1}]\times \Bbb{R}[t_1,t_2,t_1^{-1}])/\Bbb{R}1_{[0]}\eqno(4.49)$$
 with its Lie bracket induced by
$$[f_{[0]},g_{[0]}]=(t_1f_{t_1}g_{t_2}-t_1f_{t_2}g_{t_1})_{[0]},\;\;[f_{[1]},g_{[1]}]=(t_1^nfg)_{[0]},\eqno(4.50)$$
$$[f_{[0]},g_{[1]}]=[t_1f_{t_1}g_{t_2}-f_{t_2}(t_1g_{t_1}+ng/2)]_{[1]}\eqno(4.51)$$
for $f,g\in \Bbb{R}[t_1,t_2,t_1^{-1}]$ (cf. (1.23)). 
The space 
$$\bar{\cal B}=(\Bbb{R}[t_1,t_2,t_1^{-1}]\times \Bbb{R}[t_1,t_2,t_1^{-1}])/\Bbb{R}(t_1)_{[0]}\eqno(4.52)$$
with its Lie bracket induced by
$$[f_{[0]},g_{[0]}]=[f_{t_2}(t_1g_{t_1}-g)+(f-t_1f_{t_1})g_{t_2}]_{[0]},\;\;[f_{[1]},g_{[1]}]=(t_1^nfg)_{[0]},\eqno(4.53)$$
$$[f_{[0]},g_{[1]}]=[f_{t_2}(t_1g_{t_1}+(n-1)g/2)+(f-t_1f_{t_1})g_{t_2}]_{[1]}\eqno(4.54)$$
for $f,g\in \Bbb{R}[t_1,t_2,t_1^{-1}]$ (cf. (1.23)). 
\pse

{\bf Example 4.5}. Let $n$ be an integer. The space $\bar{\cal B}=\Bbb{R}[t_1,t_2,t_3,t_1^{-1}]\times \Bbb{R}[t_1,t_2,t_3,t_1^{-1}]$ with the Lie bracket:
$$[f_{[0]},g_{[0]}]=[(t_1f_{t_1}+f_{t_3})(g_{t_2}-g)+(f-f_{t_2})(t_1g_{t_1}+g_{t_3})]_{[0]},\;\;[f_{[1]},g_{[1]}]=(t_1^nfg)_{[0]},\eqno(4.55)$$
$$[f_{[0]},g_{[1]}]=[(t_1f_{t_1}+f_{t_3})(g_{t_2}-g/2)+(f-f_{t_2})(t_1g_{t_1}+g_{t_3}+ng/2)]_{[1]}\eqno(4.56)$$
for $f,g\in \Bbb{R}[t_1,t_2,t_3,t_1^{-1}]$ (cf. (1.23)).
The space 
$$\bar{\cal B}=(\Bbb{R}[t_1,t_2,t_3,t_1^{-1}]\times \Bbb{R}[t_1,t_2,t_3,t_1^{-1}])/\Bbb{R}1_{[0]}\eqno(4.57)$$
with its Lie bracket induced by
$$[f_{[0]},g_{[0]}]=[(t_1f_{t_1}+f_{t_3})g_{t_2}-f_{t_2}(t_1g_{t_1}+g_{t_3})]_{[0]},\;\;[f_{[1]},g_{[1]}]=(t_1^nfg)_{[0]},\eqno(4.58)$$
$$[f_{[0]},g_{[1]}]=[(t_1f_{t_1}+f_{t_3})g_{t_2}-f_{t_2}(t_1g_{t_1}+g_{t_3}+ng/2)]_{[1]}\eqno(4.59)$$
for $f,g\in \Bbb{R}[t_1,t_2,t_3,t_1^{-1}]$ (cf. (1.23)). 
The space 
$$\bar{\cal B}=(\Bbb{R}[t_1,t_2,t_3,t_1^{-1}]\times \Bbb{R}[t_1,t_2,t_3,t_1^{-1}])/\Bbb{R}(t_1)_{[0]}\eqno(4.60)$$
with its Lie bracket induced by
$$[f_{[0]},g_{[0]}]=[f_{t_2}(t_1g_{t_1}+g_{t_3}-g)+(f-t_1f_{t_1}-f_{t_3})g_{t_2}]_{[0]},\;\;[f_{[1]},g_{[1]}]=(t_1^nfg)_{[0]},\eqno(4.61)$$
$$[f_{[0]},g_{[1]}]=[f_{t_2}(t_1g_{t_1}+g_{t_3}+(n-1)g/2)+(f-t_1f_{t_1}-f_{t_3})g_{t_2}]_{[1]}\eqno(4.62)$$
for $f,g\in \Bbb{R}[t_1,t_2,t_1^{-1}]$ (cf. (1.23)). 
\pse

{\bf Example 4.6}. The space 
$$\bar{\cal B}=(\Bbb{R}[t_1,t_2,t_1^{-1},t_2^{-1}]\times \Bbb{R}[t_1,t_2,t_1^{-1},t_2^{-1}])/\Bbb{R}1_{[0]}\eqno(4.63)$$
with its Lie bracket induced by
$$[f_{[0]},g_{[0]}]=[t_1t_2(f_{t_1}g_{t_2}-f_{t_2}g_{t_1})]_{[0]},\;\;[f_{[1]},g_{[1]}]=(t_1t_2fg)_{[0]},\eqno(4.64)$$
$$[f_{[0]},g_{[1]}]=[t_1f_{t_1}(t_2g_{t_2}+g/2)-t_2f_{t_2}(t_1g_{t_1}+g/2)]_{[1]}\eqno(4.65)$$
for $f,g\in \Bbb{R}[t_1,t_2,t_1^{-1},t_2^{-1}]$ (cf. (1.23)).
The space 
$$\bar{\cal B}=(\Bbb{R}[t_1,t_2,t_1^{-1},t_2^{-1}]\times \Bbb{R}[t_1,t_2,t_1^{-1},t_2^{-1}])/\Bbb{R}(t_2)_{[0]}\eqno(4.66)$$
 with its Lie bracket induced by
$$[f_{[0]},g_{[0]}]=[t_1f_{t_1}(t_2g_{t_2}-g)+t_1(f-t_2f_{t_2})g_{t_1}]_{[0]},\;\;[f_{[1]},g_{[1]}]=((t_1t_2)fg)_{[0]},\eqno(4.67)$$
$$[f_{[0]},g_{[1]}]=[t_1t_2f_{t_1}g_{t_2}+(f-t_2f_{t_2})(t_1g_{t_1}+g/2)]_{[1]}\eqno(4.68)$$
for $f,g\in \Bbb{R}[t_1,t_2,t_1^{-1},t_2^{-1}]$ (cf. (1.23)).
\pse

{\bf Example 4.7}. Let $m$ and $n$ be an integer. 
The space 
$$\bar{\cal B}=(\Bbb{R}[t_1,t_2,t_3,t_1^{-1},t_2^{-1}]\times \Bbb{R}[t_1,t_2,t_3,t_1^{-1},t_2^{-1}])/\Bbb{R}1_{[0]}\eqno(4.69)$$
with its Lie bracket induced by
$$[f_{[0]},g_{[0]}]=[t_2(t_1f_{t_1}+f_{t_3})g_{t_2}-t_2f_{t_2}(t_1g_{t_1}+g_{t_3})]_{[0]},\;\;[f_{[1]},g_{[1]}]=(t_1^mt_2^nfg)_{[0]},\eqno(4.70)$$
$$[f_{[0]},g_{[1]}]=[(t_1f_{t_1}+f_{t_3})(t_2g_{t_2}+ng/2)-t_2f_{t_2}(t_1g_{t_1}+g_{t_3}+mg/2)]_{[1]}\eqno(4.71)$$
for $f,g\in \Bbb{R}[t_1,t_2,t_3,t_1^{-1},t_2^{-1}]$ (cf. (1.23)). 
The space 
$$\bar{\cal B}=(\Bbb{R}[t_1,t_2,t_3,t_1^{-1},t_2^{-1}]\times \Bbb{R}[t_1,t_2,t_3,t_1^{-1},t_2^{-1}])/\Bbb{R}(t_1)_{[0]}\eqno(4.72)$$
 with its Lie bracket induced by
$$[f_{[0]},g_{[0]}]=[t_2f_{t_2}(t_1g_{t_1}+g_{t_3}-g)+t_2(f-t_1f_{t_1}-f_{t_3})g_{t_2}]_{[0]},\;\;[f_{[1]},g_{[1]}]=(t_1^mt_2^nfg)_{[0]},\eqno(4.73)$$
$$[f_{[0]},g_{[1]}]=[t_2f_{t_2}(t_1g_{t_1}+g_{t_3}+(m-1)g/2)+(f-t_1f_{t_1}-f_{t_3})(t_2g_{t_2}+ng/2)]_{[1]}\eqno(4.74)$$
for $f,g\in \Bbb{R}[t_1,t_2,t_3,t_1^{-1},t_2^{-1}]$ (cf. (1.23)). 
The space 
$$\bar{\cal B}=(\Bbb{R}[t_1,t_2,t_3,t_1^{-1},t_2^{-1}]\times \Bbb{R}[t_1,t_2,t_3,t_1^{-1},t_2^{-1}])/\Bbb{R}(t_2)_{[0]}\eqno(4.75)$$
 with its Lie bracket induced by
$$[f_{[0]},g_{[0]}]=[(t_1f_{t_1}+f_{t_3})(t_2g_{t_2}-g)+(f-t_2f_{t_2})(t_1g_{t_1}+g_{t_3}]_{[0]},\;\;[f_{[1]},g_{[1]}]=(t_1^mt_2^nfg)_{[0]},\eqno(4.76)$$
$$[f_{[0]},g_{[1]}]=[(t_1f_{t_1}+f_{t_3})(t_2g_{t_2}+(n-1)g/2)+(f-t_2f_{t_2})(t_1g_{t_1}+g_{t_3}+mg/2)]_{[1]}\eqno(4.77)$$
for $f,g\in \Bbb{R}[t_1,t_2,t_3,t_1^{-1},t_2^{-1}]$ (cf. (1.23)). 
\pse

{\bf Example 4.8}. Let $k$ be a positive integer and let ${\cal A}$ be the subalgebra of 
$$\Bbb{R}[t_1^{1/k},t_2^{1/k},t_3,t_4, t_1^{-1/k},t_2^{-1/k}]\eqno(4.78)$$
 generated by 
$$\{t_1,t_2,t_3,t_4,(t_1t_2)^{1/k},t_1^{-1},t_2^{-1}\}.\eqno(4.79)$$
In particular, ${\cal A}=\Bbb{R}[t_1,t_2,t_3,t_4, t_1^{-1},t_2^{-1}]$ when $k=1$. Suppose that $m$ and $n$ are two fixed integers. The space 
$$\bar{\cal B}=({\cal A}\times {\cal A})/\Bbb{R}1_{[0]}\eqno(4.80)$$
 with its Lie bracket induced by
$$[f_{[0]},g_{[0]}]=[(t_1f_1+f_{t_3})(t_2g_{t_2}+g_{t_4})-(t_2f_{t_2}+f_{t_4})(t_1g_{t_1}+g_{t_3})]_{[0]},\eqno(4.81)$$
$$[f_{[1]},g_{[1]}]=(t_1^mt_2^nfg)_{[0]},\eqno(4.82)$$
$$[f_{[0]},g_{[1]}]=[(t_1f_{t_1}+f_{t_3})(t_2g_{t_2}+g_{t_4}+ng/2)-(t_2f_{t_2}+f_{t_4})(t_1g_{t_1}+mg/2)]_{[1]}\eqno(4.83)$$
for $f,g\in {\cal A}$ (cf. (1.23)).
The space 
$$\bar{\cal B}=({\cal A}\times {\cal A})/\Bbb{R}(t_2)_{[0]}\eqno(4.84)$$
with its Lie bracket induced by
$$[f_{[0]},g_{[0]}]=[(t_1f_{t_1}+f_{t_3})(t_2g_{t_2}+g_{t_4}-g)+(f-t_2f_{t_2}-f_{t_4})(t_1g_{t_1}+g_{t_3})]_{[0]},\eqno(4.85)$$
$$[f_{[1]},g_{[1]}]=((t_1^mt_2^n)fg)_{[0]},\eqno(4.86)$$
\begin{eqnarray*}[f_{[0]},g_{[1]}]&=&[(t_1f_{t_1}+f_{t_3})(t_2g_{t_2}+g_{t_4}+(n-1)g/2)+(f-t_2f_{t_2}-f_{t_4})\\& &\times (t_1g_{t_1}+g_{t_3}+mg/2)]_{[1]}\hspace{8.1cm}(4.87)\end{eqnarray*}
for $f,g\in{\cal A}$ (cf. (1.23)).

\vspace{0.4cm}

\begin{center}{\Large \bf Acknowledgments}\end{center}
\vspace{0.3cm}

We thank the referee for his/her pointing out some misprints and making suggestions of corrections.

\vspace{1cm}

\noindent{\Large \bf References}

\hspace{0.5cm}

\begin{description}

\item[{[B]}] R. Block, On torsion-free abelian groups and Lie algebras, {\it Proc. Amer. Math. Soc.} {\bf 9} (1958), 613-620.

\item[{[DZ]}] D. Dokovic and K. Zhao, Derivations, isomorphisms and second cohomology of generalized Block algebras, {\it Algebra Colloq.} {\bf 3} (1996), 245-272.

\item[{[X1]}] X, Xu, On simple Novikov algebras and their irreducible modules, {\it J. Algebra} {\bf 185} (1996), 905-934.

\item[{[X2]}] X. Xu, Novikov-Poisson algebras, {\it J. Algebra} {\bf 190} (1997), 253-279.

\item[{[X3]}] X. Xu, Quadratic conformal superalgebras, {\it submitted}.

\item[{[X4]}] X. Xu, New simple Lie algebras of generalized Cartan types over a field with characteristic 0, {\it J. Algebra}, to appear.

\item[{[X5]}] X. Xu, {\it Introduction to Vertex Operator Superalgebras and Their Modules}, Kluwer Academic Publishers, 1998.

\item[{[X6]}] X. Xu, {\it Algebraic Theory of Hamiltonian Superoperators}, to appear.

\end{description}
\end{document}